\documentclass[]{siamltex1213}
\pdfoutput=1

\usepackage{amsmath}
\usepackage{amsfonts}
\usepackage{graphicx}
\usepackage{algorithmic}

\usepackage{amssymb}

\newcommand{\TheTitle}{Degenerate Kalman filter error covariances and their convergence onto the unstable subspace\thanks{This paper has been published by the SIAM/ASA Journal of Uncertainty Quantification (Vol. 5, Issue 1, 2017, pp. 304--333); \url{http://epubs.siam.org/toc/sjuqa3/5/1}; doi: 10.1137/16M1068712.} }
\newcommand{\TheAuthors}{M. Bocquet, K. S. Gurumoorthy, A. Apte, A. Carrassi, C. Grudzien, and C. K. R. T. Jones}

\title{{\TheTitle}}

\author{
  Marc Bocquet\thanks{CEREA, joint laboratory \'Ecole des Ponts ParisTech and EDF R\&D, Universit\'e Paris-Est, Champs-sur-Marne, France (\email{marc.bocquet@enpc.fr})}
  \and
  Karthik S. Gurumoorthy\thanks{International Center for Theoretical Sciences, Tata Institute of Fundamental Research, Bangalore, Karnataka, India}
  \and
  Amit Apte\footnotemark[3]
  \and
  Alberto Carrassi\thanks{Nansen Environmental and Remote Sensing Center, Bergen, Norway}
  \and
  Colin Grudzien\thanks{Department of Mathematics, University of North Carolina, Chapel Hill, North Carolina, USA. Current address: Nansen Environmental and Remote Sensing Center, Bergen, Norway} 
  \and
  Christopher K. R. T. Jones\thanks{Department of Mathematics, University of North Carolina, Chapel Hill, North Carolina, USA.} 
}

\usepackage{amsopn}

\ifpdf
\hypersetup{
  pdftitle={\TheTitle},
  pdfauthor={\TheAuthors}
}
\fi

\def\be{\begin{equation}}
\def\ee{\end{equation}}
\def\bea{\begin{eqnarray}}
\def\eea{\end{eqnarray}}
\def\nn{\nonumber\\}
\def\({\left(}
\def\){\right)}

\def\bx{{\mathbf x}}
\def\by{{\mathbf y}}

\def\bw{{\mathbf w}}
\def\bu{{\mathbf u}}
\def\bv{{\mathbf v}}

\def\bh{{\mathbf h}}

\def\bA{\mathbf A}
\def\bB{\mathbf B}
\def\bC{\mathbf C}
\def\wbC{\widetilde{\mathbf C}}
\def\bD{\mathbf D}
\def\bE{\mathbf E}

\def\bL{\mathbf L}
\def\bG{\mathbf G}

\def\bR{\mathbf R}
\def\bS{\mathbf S}
\def\bP{\mathbf P}
\def\bM{\mathbf M}
\def\bH{\mathbf H}
\def\bI{\mathbf I}
\def\bQ{\mathbf Q}
\def\bT{\mathbf T}

\def\bU{\mathbf U}
\def\bV{\mathbf V}
\def\bX{\mathbf X}
\def\bY{\mathbf Y}
\def\bW{\mathbf W}
\def\bZ{\mathbf Z}

\def\bJ{\mathbf J}

\def\bOmega{{\boldsymbol \Omega}}
\def\bomega{{\boldsymbol \omega}}
\def\bGamma{{\boldsymbol \Gamma}}
\def\bPsi{{\boldsymbol \Psi}}
\def\bPhi{{\boldsymbol \Phi}}
\def\bTheta{{\boldsymbol \Theta}}
\def\bSigma{{\boldsymbol \Sigma}}
\def\bXi{{\boldsymbol \Xi}}

\def\bLambda{{\boldsymbol \Lambda}}

\def\bzero{{\mathbf 0}}
\def\b1{{\mathbf 1}}

\def\T{{\mathrm T}}
\def\ve{(\varepsilon)}
\def\Rn{{\mathbb R}^n}
\newcommand\cone[1]{{\mathcal C}^{#1}}
\newcommand\conep[1]{{\mathcal C}^{#1}_{+}}
\newcommand\R[2]{{\mathbb R}^{#1 \times #2}}

\def\Ul{{\mathcal U}_k}
\def\Sl{{\mathcal S}_k}

\def\Sv{{\mathcal S}_{k:0}}
\def\Svi{{\mathfrak S}_k}
\def\Sli{{\mathfrak S}}
\def\rank{\mathrm{rank}}
\def\im{\mathrm{Im}}
\def\W{{\mathcal W}}
\def\V{{\mathcal V}}
\def\s{{s}}

\def\Gpp{\bG_{\!_{++}}}
\def\Gpm{\bG_{\!_{+-}}}
\def\Gmp{\bG_{\!_{-+}}}
\def\Gmm{\bG_{\!_{--}}}

\def\Tpp{\bTheta_k^{\!^{++}}}
\def\Tpm{\bTheta_k^{\!^{+-}}}
\def\Tmp{\bTheta_k^{\!^{-+}}}
\def\Tmm{\bTheta_k^{\!^{--}}}

\def\Gapp{\bGamma_{\!^{++},k}}

\newcommand\La[1]{\bLambda_{#1}}
\newcommand\Lp[1]{\bLambda_{+,#1}}
\newcommand\Lm[1]{\bLambda_{-,#1}}
\newcommand\Lpp[1]{\bLambda_{\!_{++},#1}}
\newcommand\Lpn[1]{\bLambda_{\!_{+0},#1}}
\newcommand\Lpg[1]{\bLambda_{\!_{+\uparrow},#1}}
\newcommand\Lpl[1]{\bLambda_{\!_{+\downarrow},#1}}
\newcommand\Lpng[1]{\bLambda_{\!_{+0\uparrow},#1}}
\newcommand\Lpnl[1]{\bLambda_{\!_{+0\downarrow},#1}}

\newcommand\I[1]{\bI_{#1}}

\newcommand\mat[2]{\left[ \begin{array}{cc} #1 & \bzero \\ \bzero & #2 \end{array} \right]}

\usepackage{xcolor}

\begin{document}

\maketitle

\begin{abstract}
The characteristics of the model dynamics are critical in the performance of (ensemble) Kalman filters. In particular,
as emphasized in the seminal work of Anna Trevisan and co-authors, the error covariance matrix is asymptotically
supported by the unstable-neutral subspace only, i.e., it is spanned by the backward Lyapunov vectors with
non-negative exponents. This behavior is at the core of algorithms known as Assimilation in the Unstable Subspace,
although a formal proof was still missing.

This paper provides the analytical proof of the convergence of the Kalman filter covariance matrix onto the
unstable-neutral subspace when the dynamics and the observation operator are linear and when the dynamical model is
error-free, for any, possibly rank-deficient, initial error covariance matrix. The rate of convergence is provided as
well.  The derivation is based on an expression that explicitly relates the error covariances at an arbitrary time to
the initial ones.  It is also shown that if the unstable and neutral directions of the model are sufficiently observed
and if the column space of the initial covariance matrix has a non-zero projection onto all of the forward Lyapunov
vectors associated with the unstable and neutral directions of the dynamics, the covariance matrix of the Kalman filter
collapses onto an asymptotic sequence which is independent of the initial covariances.  Numerical results are also shown
to illustrate and support the theoretical findings.
\end{abstract}



\section{Introduction}
\label{sec:intro}

\subsection{Context and objectives}

Filtering methods are the techniques of estimation theory that process measurements sequentially as they become
available.  In a probabilistic Bayesian framework, the output of a filter is a probability density function (pdf), the
conditional posterior pdf $p(\bx\vert\by)$ of the process $\bx$, given the data $\by$ and a prior distribution
$p(\bx)$.  The posterior pdf fully characterizes the state's estimation and quantifies the uncertainty of the
estimate. However, its exact calculation is extremely difficult in practice, and most often computationally intractable
in high-dimensional, complex systems, such as numerical climate and weather models.

For linear dynamics, measurements with a linear dependence on the state variables, and Gaussian errors, the Kalman
filter (KF) is the optimal filtering solution \cite{kalman1960}.  The Gaussian hypothesis implies an enormous
simplification: the pdfs are all completely characterized by their first and second moments. In this case, the error
covariance matrix quantifies the uncertainty of the state's estimate represented by the mean.  The KF has been extremely
successful for decades in numerous fields including navigation, economy, robotics, tracking objects, adaptive optics, and
many computer vision applications.

A Monte Carlo formulation of the KF leads to the introduction of a class of {\it Gaussian} algorithms referred to as
ensemble Kalman filters (EnKFs) \cite{evensen2009}. They have been widely applied in atmospheric and oceanic contexts,
where all methods designed for filtering or smoothing are referred to as data assimilation (DA).  In the EnKF the
transition probability of the process, as well as all the error covariances entering the assimilation of observations,
are approximated using an ensemble of realizations (members in the EnKF jargon) of the model dynamics. The EnKF and its
variants are currently among the most popular approaches for DA in high-dimensional systems. Evidence has emerged that a
small number of members, typically 100, is sufficient in many applications, especially when using localization
techniques \cite[and references therein]{sakov2011}, hence making the EnKF feasible in situations where the forward step
of DA is computationally expensive. The choice of the ensemble members is critical and a key aspect in the EnKF
setup. While a large ensemble is generally desirable to explain and represent the actual uncertainty in the most
realistic manner, their number is limited by the computational resources at disposal. In the absence of localization,
the EnKF error covariances are thus degenerate (or rank-deficient) by construction and it is then relevant to adequately
choose these few (much fewer than the system's dimension) members so as to maximize the representation of the actual
uncertainty.

For nonlinear chaotic dynamics, the assimilation in the unstable subspace (AUS), introduced by Anna Trevisan and
collaborators \cite{trevisan2004, carrassi2008a, trevisan2010, trevisan2011, palatella2015}, has shed light on an
efficient way to operate the assimilation of observations by using the unstable subspace to describe the uncertainty in
the estimate.  AUS is based on two key properties of deterministic, typically dissipative, chaotic systems: (i) the
perturbations tend to project on the unstable manifold of the dynamics, and (ii) the dimension of the unstable manifold
is typically much smaller than the full phase-space dimension.  Applications to atmospheric, oceanic, and traffic models
\cite{carrassi2008b, uboldi2006,palatella2013b} showed that even in high-dimensional systems, an efficient error control
is achieved by monitoring only the unstable directions, and sometimes even a subset of them, making AUS a
computationally efficient alternative to standard procedures.

The AUS approach has recently motivated a research effort toward a proper mathematical formulation and assessment of its
driving idea, i.e., the span of the estimation error covariance matrices asymptotically (in time) tends to the
subspace spanned by the unstable and neutral backward Lyapunov vectors (BLVs). A proper statement of this latter property in
precise mathematical terms is of vast importance for the design of efficient reduced-order uncertainty quantification
and DA methods.

The first recent result along this line is given in \cite{gurumoorthy2017}. It is proved that for linear, discrete,
autonomous and non-autonomous, deterministic systems (perfect model) with noisy observations, the covariance equations in
the KF asymptotically bound the rank of the forecast and the analysis error covariance matrices to be less than or
equal to the number of non-negative Lyapunov exponents of the system. Further, the support of these error covariance
matrices is shown to be confined to the space spanned by the unstable and neutral BLVs.  The
results in \cite{gurumoorthy2017} were obtained assuming a full rank covariance matrix at initial time. The conditions
that imply the convergence, for possibly degenerate (low rank) initial matrices remained unresolved, yet they are
fundamental to link these mathematical findings with concrete reduced-rank DA methods, particularly the EnKF.

This is the subject of the present work, which studies the convergence in the general setting of degenerate covariance
matrices.  A pivotal result is the analytic proof of the KF covariance collapse, for any initial error covariance (of
arbitrary rank), onto the unstable-neutral subspace. We also provide rigorous mathematical results for the rate of
convergence on the stable subspace and for the asymptotic rank of the error covariance matrix.  Finally, we derive an
expression for the asymptotic error covariance matrix as a function of the initial one. This in turn allows us to prove,
under certain observability conditions, the existence of an asymptotic sequence of error covariance matrices, which is
independent of the initial condition.

In the following, we set up the notations and discuss the organization of the paper.

\subsection{Problem formulation}
\label{sec:theory}

The purpose of this paper is the estimation of the unknown state of a system based on partial and noisy observations. The dynamical
and observational models are both assumed to be linear, and expressible as
\begin{align}
  \label{eq:dynmodel}
  \bx_{k} &= \bM_{k}\bx_{k-1} + \bw_k , \\
  \label{eq:obsmodel}
  \by_{k} &= \bH_k \bx_k + \bv_k ,
 \end{align}
with $\bx\in{\mathbb R}^n$ and $\by\in{\mathbb R}^d$ being the system's state and observation, respectively, related via
the linear observation operator $\bH_k: {\mathbb R}^n \mapsto {\mathbb R}^d$.  Throughout the entire text the
conventional notation $k=0,1,2,\ldots$ is adopted to indicate that the quantity is defined at time $t_k$.  The matrix
$\bM_{k:l}$ is taken to represent the resolvent of the linear forward model from time $t_{l}$ to time $t_k$, and is
assumed to be non-singular throughout this paper. In particular $\bM_{k:k}=\I{n}$, where $\I{n}$ is the identity matrix (of
size $n\times n$ in this case). We designate $\bM_k$ as the one-step matrix resolvent of the forward model from
$t_{k-1}$ to $t_k$: $\bM_k \triangleq \bM_{k:k-1}$ and, consequently, $\bM_{k:l} = \bM_{k}\bM_{k-1} \ldots \bM_{l+1}$, with the symbol
$\triangleq$ used to signify that the expression is a definition.  We will assume that the Lyapunov spectrum of the
dynamics defined by $\bM_{k:0}$ is non-degenerate, i.e., the Lyapunov exponents are all distinct.  This assumption
substantially simplifies the derivations that follow.  Nonetheless, most of the results in this paper can be generalized
to the degenerate case.

The model and observation noise, $\bw_k$ and $\bv_k$, are assumed mutually independent, unbiased Gaussian white
sequences with statistics
\be
   {\rm E}[\bv_k\bv_l^\T] = \delta_{k,l}\bR_k , \quad  {\rm E}[\bw_k\bw_l^\T] = \delta_{k,l}\bQ_k , \quad {\rm E}[\bv_k\bw_l^\T] = \bzero ,
\ee
where $\bR_k\in {\mathbb R}^{d\times d}$ is the observation error covariance matrix at time $t_k$, and $\bQ_k \in
{\mathbb R}^{n\times n}$ stands for the model error covariance matrix. The error covariance matrix $\bR_k$ can be
assumed invertible without losing generality.

The forecast error covariance matrix $\bP_k$ of the KF satisfies the following recurrence, the discrete-time
dynamic Riccati equation \cite{kalman1960,gurumoorthy2017}
\be
\label{eq:recurrence}
  \bP_{k+1} = \bM_{k+1}\(\I{n}+\bP_k \bOmega_k\)^{-1}\bP_k\bM_{k+1}^{\T} +\bQ_{k+1} ,
\ee
where
\be
\bOmega_k \triangleq \bH_k^{\T}\bR_k^{-1}\bH_k
\ee
is the precision matrix of the observations transferred in state space.  To avoid pathological behaviors, we will assume
in this paper that the $\left\{ \bOmega_k \right\}_{k=0,1,\ldots}$ are uniformly bounded from above, which is a very
mild hypothesis.

Equation~\eqref{eq:recurrence} highlights that the error covariance matrix, $\bP_{k+1}$, is the result of a two-step
process, consisting of the \emph{update} or \emph{analysis} step at time $t_k$ leading to the analysis error covariance
matrix $\bP^{\rm a}_k$,
\be
\label{eq:kfupdate} 
\bP^{\rm a}_k = \(\I{n}+\bP_k \bOmega_k\)^{-1}\bP_k ,
\ee
and the \emph{forecast} step which consists of
the forward propagation of the analysis error covariance,
\be
\label{eq:kfforecast} 
\bP_{k+1} = \bM_{k+1}\bP^{\rm a}_k\bM_{k+1}^{\T} + \bQ_{k+1} .
\ee
It is worth mentioning that \eqref{eq:recurrence} still holds when $\bP_k$ is degenerate, i.e., $\rank(\bP_k) <
n$.  This is the typical circumstance encountered in the EnKF \cite[and references
  therein]{evensen2009}.  In this case, assuming that the model is perfect ($\bQ_k=\bzero$) and under the same
hypotheses of linear observation and evolution operators as well as of Gaussian statistics for the initial condition and
observational errors, \eqref{eq:recurrence} will apply to the EnKF too.

\subsection{Outline of the paper}

In the rest of the paper, we will refer to \eqref{eq:recurrence} as the recurrence equation for $\bP_k$, although we
will mostly study the {\em perfect dynamical model case}, in which $\bQ_k=\bzero$.  In section \ref{sec:PkP0} we
demonstrate a relation between $\bP_k$ at any arbitrary time, $t_k>t_0$, and the initial error covariance matrix
$\bP_0$, in the general case with $\bP_0$ possibly being degenerate. An alternative proof based on the linear symplectic
representation of the KF is proposed in appendix \ref{app:symplectic}.  In the following section \ref{sec:bound}, we
derive a useful bound that plays a central role in all results and derivations discussed in this study. Then in section
\ref{sec:PkCONS} we study the asymptotic behavior of $\bP_k$ (for $k \rightarrow \infty$) along with other relevant
properties.  Section \ref{sec:asymptote} provides a proof, using a condition on the initial $\bP_0$ and certain
observability conditions, that the error covariances collapse onto an asymptotic sequence which is independent of the
initial covariance matrix $\bP_0$. Section \ref{sec:numeric} describes the numerical results corroborating and
illustrating the theoretical findings while the conclusions are drawn in section \ref{sec:concl}.

\section{Computation of the forecast error covariance matrix $\bP_k$}
\label{sec:PkP0}

In this section we consider the perfect model case, i.e., $\bQ_k=\bzero$ for all $k$. The stochastic model
case, $\bQ_k \ne \bzero$, is briefly considered in section \ref{sec:bound}.

The recurrence \eqref{eq:recurrence} is rational in $\bP_k$.  Furthermore if we assume that the $\bP_k$ are
invertible, we can take the inverse of both sides of the recurrence and obtain
\be
\label{eq:information}
\bP^{-1}_{k+1} = \bM_{k+1}^{-{\T}}\( \bP^{-1}_k + \bOmega_k \)\bM_{k+1}^{-1} ,
\ee
which shows that $\bP^{-1}_{k+1}$ is an affine function of $\bP^{-1}_k$.  This relation is usually called the {\it
  information filter} \cite[section 3.2]{simon2006}.

However, a relevant situation in applications is when the $\bP_k$ are degenerate. In this case the inverse of both sides
of \eqref{eq:information} are undefined, and a suitable generalization of \eqref{eq:information} is
required. To that end, we introduce an analytic continuation of \eqref{eq:recurrence}.
A regularized $\bP_0$ is defined as
\be
\label{eq:P0def}
\bP_0\ve \triangleq \bP_0+\varepsilon\I{n}
\ee
with $\varepsilon>0$ and we define the subsequent $\bP_k\ve$ via the recurrence
\be
\label{eq:recurrence3}
\bP_{k+1}\ve \triangleq \bM_{k+1}\(\I{n}+\bP_k\ve \bOmega_k\)^{-1}\bP_k\ve\bM_{k+1}^{\T} .
\ee
From (\ref{eq:P0def},\ref{eq:recurrence3}), $\bP_k\ve$ is seen to be full-rank.
Moreover, taking the limit $\varepsilon \rightarrow 0^+$, leads $\bP_0\ve$ continuously back to $\bP_0$ and
\eqref{eq:recurrence3} to \eqref{eq:recurrence}, so that we have
\be
\lim_{\varepsilon \rightarrow 0^+} \bP_k\ve = \bP_k(0) = \bP_k .
\ee

Then, we take the inverse of both sides of \eqref{eq:recurrence3},
\begin{align}
  \bP^{-1}_{k+1}\ve &= \bM_{k+1}^{-{\T}}\bP^{-1}_k\ve \(\I{n} + \bP_k\ve\bOmega_k \)\bM_{k+1}^{-1} \nn
  & = \bM_{k+1}^{-{\T}}\bP^{-1}_k\ve\bM_{k+1}^{-1} + \bM_{k+1}^{-{\T}}\bOmega_k \bM_{k+1}^{-1} .
\end{align}
This recurrence can easily be solved and it yields
\be
\label{eq:precision}
  \bP_k^{-1}\ve = \bM_{k:0}^{-{\T}}\bP_0^{-1}\ve\bM_{k:0}^{-1} + \bGamma_k ,
\ee
where 
\be
\label{eq:gammadef}
\bGamma_k \triangleq \sum_{l=0}^{k-1} \bM_{k:l}^{-\T} \bOmega_l \bM_{k:l}^{-1} .
\ee
This matrix, known as the information matrix \cite{jazwinski1970}, is a measure of the observability of the system since
it propagates the precision matrices $\bOmega_l$ up to $t_k$, and \eqref{eq:precision} states that the precision in
the state estimate is the sum of the forecast precision in the initial condition plus the precision of the observations
transferred into the model space.

Let us now recall the partial order defined in the cone $\cone{n}$ of the symmetric positive semi-definite matrices of
$\R{n}{n}$, of which we will make great use in this study. Similarly the partial order acts in the cone $\conep{n}$ of the
symmetric positive definite matrices of $\R{n}{n}$.  We will refer to this partial order using the standard comparison
symbols. In appendix \ref{app:cone}, its definition is provided along with some additional properties that we rely on in
this study.  From \eqref{eq:precision} and using this partial order, we have
\be
\bP_k^{-1}\ve \ge \bGamma_k .
\ee
Let us assume that the system is observable, a condition defined here as $\mathrm{det}(\bGamma_k) \neq 0$ according to
\cite{jazwinski1970} and references therein.  This yields $\bP_k\ve \le \bGamma_k^{-1}$ (see appendix \ref{app:cone}, point
3), which implies
\be
\label{eq:bound2}
  \bP_k \le \bGamma_k^{-1}.
\ee

By taking the inverse of both sides of \eqref{eq:precision} we have
\begin{align}
  \bP_k\ve &= \( \bM_{k:0}^{-{\T}}\bP_0^{-1}\ve\bM_{k:0}^{-1} + \bGamma_k \)^{-1} \nn
& = \left[ \(\I{n} + \bGamma_k \bM_{k:0}\bP_0\ve\bM_{k:0}^\T\) \bM_{k:0}^{-{\T}}\bP_0^{-1}\ve\bM_{k:0}^{-1} \right]^{-1} \nn
& = \bM_{k:0}\bP_0\ve\bM_{k:0}^\T \(\I{n} + \bGamma_k \bM_{k:0}\bP_0\ve\bM_{k:0}^\T\)^{-1}.
\end{align}
The limit $\varepsilon \rightarrow 0^+$ finally leads to
\be
\label{eq:final2}
  \bP_k = \bM_{k:0}\bP_0\bM_{k:0}^\T \(\I{n} + \bGamma_k \bM_{k:0}\bP_0\bM_{k:0}^\T\)^{-1} .
\ee
Equation~\eqref{eq:final2} is extremely important as it directly relates $\bP_k$ to $\bP_0$.  In particular it shows
that $\bP_k$ depends on two concurring factors, the matrix $\bGamma_k$ encoding all information about the observability
of the system, and the matrix $\bM_{k:0}\bP_0 \bM_{k:0}^{\T}$ representing the free forecast of the initial
covariances. The latter exemplifies the uncertainty propagation under the model dynamics, the former the ability of the
observations to counteract the error growth.

We now use the matrix shift lemma that asserts that for any matrices $\bA \in \R{l}{m}$ and $\bB \in \R{m}{l}$, we have
$\bA f(\bB\bA)=f(\bA\bB)\bA$, with $x \mapsto f(x)$ being any function that can be expressed as a formal power series.
A derivation is recalled in appendix \ref{app:sml}.  Here, we choose $f(x)=(1+x)^{-1}$, $\bA=\bM_{k:0}^\T$ and $\bB=\bGamma_k
\bM_{k:0}\bP_0$, to obtain an alternative formulation of \eqref{eq:final2}
\be
\label{eq:final3}
  \bP_k = \bM_{k:0}\bP_0  \left[ \I{n} + \bM_{k:0}^\T\bGamma_k \bM_{k:0}\bP_0 \right]^{-1}\bM_{k:0}^{\T}
\ee
or, in a more condensed form,
\be
\label{eq:final}
  \bP_k = \bM_{k:0}\bP_0  \left[ \I{n} + \bTheta_k\bP_0 \right]^{-1}\bM_{k:0}^{\T}
\ee
where 
\be
\label{eq:thetadef}
\bTheta_k \triangleq \bM_{k:0}^{\T} \bGamma_k \bM_{k:0} = \sum_{l=0}^{k-1} \bM_{l:0}^{\T} \bOmega_l \bM_{l:0} .
\ee
This matrix is also related to the observability of the system
but pulled back at the initial time $t_0$.

A more general, albeit less straightforward, proof of the expressions for $\bP_k$ as a function of $\bP_0$ can be
obtained using the underlying symplectic structure of the KF and is described in appendix \ref{app:symplectic}.

\section{Free forecast of $\bP_0$ as an upper bound}
\label{sec:bound}

In this section we demonstrate that an upper bound for $\bP_k$ is given by the {\em free} forecast of $\bP_0$; the term
free is used in this study to mean without the observational forcing applied at analysis times.  It is worth mentioning
already that, although the existence of this bound is indeed very intuitive, its formal proof is provided here because
it plays a pivotal role in all the convergence results that follow.  The bound can be derived directly from the general
expression for $\bP_k$, \eqref{eq:final}, but we opted for showing a different approach, independent from
\eqref{eq:final}, that better highlights the relevance of the bound for the results that follow.

The error covariance matrix $\bP_k$ is symmetric and our purpose is to make the recurrence equation look patently
symmetric as well so that we can derive inequalities using the partial ordering in $\cone{n}$.  As a positive
semi-definite matrix, $\bP_k$ can be decomposed into $\bP_k=\bX_k\bX_k^\T$ using, for instance, a Choleski
decomposition, with $\bX_k\in\R{n}{m}$ ($m\le n$). Here, as opposed to the rest of the paper and for the sake of
generality, we consider the presence of model noise given that it only represents a minor complication.  The recurrence
equation can be written as
\be
\bP_{k+1} = \bM_{k+1}\(\I{n}+\bX_k\bX_k^\T \bOmega_k\)^{-1}\bX_k\bX_k^\T\bM_{k+1}^{\T} + \bQ_{k+1}.
\ee
We use again the matrix shift lemma but this time with $f(x)=(1+x)^{-1}$, $\bA =
\bX_k$, and $\bB=\bX_k^\T \bOmega_k$ so that \eqref{eq:recurrence} becomes
\be
\label{eq:recurrence2}
  \bP_{k+1} = \bM_{k+1}\bX_k\(\I{m}+\bX_k^\T \bOmega_k\bX_k\)^{-1}\bX_k^\T\bM_{k+1}^{\T} + \bQ_{k+1}.
\ee
Using the partial order in $\cone{m}$, we have from
\be
\(\I{m}+\bX_k^\T \bOmega_k\bX_k\)^{-1} \le \I{m}
\ee
and from \eqref{eq:recurrence2} that
\be
\bQ_{k+1} \le \bP_{k+1} \le \bM_{k+1}\bP_k\bM_{k+1}^{\T} + \bQ_{k+1}.
\ee
Hence $\bP_k$ is bounded from above by the free forecast $\widetilde{\bP}_k$
that satisfies $\widetilde{\bP}_0 = \bP_0$ and the recurrence
\be
\widetilde{\bP}_{k+1} = \bM_{k+1}\widetilde{\bP}_k\bM_{k+1}^{\T} + \bQ_{k+1}
\ee
whose solution is, for $k \ge 0$,
\be
\widetilde{\bP}_{k} = \bM_{k:0}\bP_0 \bM_{k:0}^{\T} + \bXi_k,
\ee
where
\be
\bXi_0 \triangleq \bzero \, \quad \text{and, for} \, k\ge 1, \quad \bXi_k \triangleq
\sum_{l=1}^k \bM_{k:l}\bQ_l \bM_{k:l}^{\T}
\ee
is known as the controllability matrix \cite{jazwinski1970}.
Therefore
\be
\label{eq:bound0}
\bQ_k \le \bP_k \le \bM_{k:0}\bP_0 \bM_{k:0}^{\T} + \bXi_k.
\ee
In particular, in the perfect model case, we obtain the pivotal inequality
\be
\label{eq:bound1}
  \bP_k \le \bM_{k:0}\bP_0 \bM_{k:0}^{\T} .
\ee
  
Under the aforementioned assumptions on linear dynamical and observational models and Gaussian error statistics, the
inequalities \eqref{eq:bound0} and \eqref{eq:bound1} state that DA will always reduce and, in the worst case, leave
unchanged, the state's estimate uncertainty with respect to the free run.  This effect was already discussed in the
context of nonlinear dynamics and in relation to the stability properties of DA systems in \cite{carrassi2008a},
although an analytic proof in the nonlinear case is not provided either in that work or in the present one.

\section{Convergence of the error covariance matrix: theoretical results}
\label{sec:PkCONS}

This section describes some of the implications of the recurrence equation and bounds described in the previous section
that are relevant for the design of reduced-order formulations of the Kalman filter with unstable dynamics. We will
assume here, again, to be in the perfect model scenario, $\bQ_k= \bzero$.

\subsection{Rank of $\bP_k$}
\label{sec:RANKP}

From the inequality \eqref{eq:bound1}, it is clear that the column space of $\bP_k$, i.e., the subspace
$\im(\bP_k) = \left\{ \bP_k\bx, \, \bx \in \Rn \right\}$ satisfies
\be
\im(\bP_k) \subseteq \bM_{k:0}\( \im(\bP_0) \) .
\ee
Moreover, since from \eqref{eq:recurrence} (with $\bQ_k=0$)
\be
\rank(\bP_{k+1}) = \rank(\bP_{k}) ,
\ee
we infer that
\be
  \im(\bP_k) = \bM_{k:0}\(\im(\bP_0)\) .
\ee
The moral is that the KF merely operates within the subspaces of the sequence
$\bM_{k:0}\(\im(\bP_0)\)$, which does not depend on the observations.
In the absence of model error, the rank of $\bP_k$ cannot exceed that of $\bP_0$ even if the dynamics are degenerate.

\subsection{Collapse of the error covariance matrices onto the unstable-neutral subspace}
\label{sec:collapse}

The unstable-neutral subspace is defined as the subspace $\Ul$ spanned by the $n_0$ BLVs at $t_k$ whose exponents,
$\lambda_i$ with $i=1,...,n_0$, are non-negative. The stable subspace $\Sl$ is defined as the subspace spanned by the
$n-n_0$ BLVs at $t_k$ associated with negative exponents.  The inequality \eqref{eq:bound1}, $\bP_k \le
\bM_{k:0}\bP_0 \bM_{k:0}^{\T}$, provides the convergence onto $\Ul$ in a sense that is made clear below. It also gives
the rate of such convergence as shown in section \ref{sec:rate}.

Let us write the singular value decomposition (SVD) of $\bM_{k:0}=\bU_{k:0}\bSigma_{k:0}\bV_{k:0}^\T$, where $\bU_{k:0}$ and
$\bV_{k:0}$ are both orthogonal matrices in $\R{n}{n}$, and $\bSigma_{k:0}$ in $\conep{n}$ is the diagonal matrix of the
singular values.  The left singular vectors are the columns of $\bU_{k:0} = [\bu^{k:0}_1, \ldots, \bu^{k:0}_n]$ and when
$k \rightarrow \infty$, they converge to the BLVs defined at $t_k$, denoted here as $\bu^k_i$. The
right singular vectors are the columns of $\bV_{k:0} = [\bv^{k:0}_1, \ldots, \bv^{k:0}_n]$ which converge to the forward
Lyapunov vectors (FLVs) at time $t_0$ as $k \rightarrow \infty$ denoted here as $\bv^0_i$ \cite{legras1996, trevisan1998}. Let
us write $[\bSigma_{k:0}]_{i,i} = \exp(\lambda_i^k k)$ with $\lambda_i^k$ being real numbers and for large $k$ ordered
as $\lambda_1^k >\cdots > \lambda_{n_0}^k \geq 0> \lambda^k_{n_0+1}> \cdots>\lambda_n^k$, which is justified by the
non-degeneracy hypothesis on the Lyapunov spectrum.  Using the SVD we have
\be
\label{eq:decomposition}
\bM_{k:0}\bP_0 \bM_{k:0}^{\T} = \bU_{k:0}\bSigma_{k:0}\bV_{k:0}^\T\bP_0\bV_{k:0}\bSigma_{k:0}\bU^\T_{k:0} .
\ee

Define $\Svi$ as the set of indices $i$ for which $\lambda^k_i < 0$ and $\Svi^s$ to be the set of indices corresponding
to the $\s$ smallest singular values in $\bSigma_{k:0}$. Note also that $\Svi^{n-n_0} = \Svi$ for large $k$. Let the
subspace $\Sv^s$ be the span of the left singular vectors $\bu^{k:0}_i$, where $i \in \Svi^s$. Let $\Pi_{\Sv^s}$ be the
orthogonal projector onto $\Sv^s$ which, owing to the orthonormality of the left singular vectors, reads
\be
\Pi_{\Sv^s} = \sum_{i \in \Svi^s} \bu_i^{k:0}\(\bu_i^{k:0}\)^{\T} .
\ee
For large enough $k$, $\Svi$ gets progressively closer and eventually coincides with the set $\Sli$ of indices $i$ for
which $\lambda_i < 0$ and each of its subsets $\Svi^s$ approaches its corresponding subset $\Sli^s$ defined
similarly. Note that $\Sli^{n-n_0} = \Sli$. Furthermore, the subspace ${\Sv^s}$ converges to $\Sl^s$ which is the span
of the $\s$ most stable BLVs.

We are now interested in an upper bound in $\cone{s}$ for $\(\bV^s_{k:0}\)^\T\bP_0\bV^s_{k:0}$, with $\bV^s_{k:0} =
[\bv^{k:0}_{n-s+1}, \ldots, \bv^{k:0}_n]$ to be jointly used with \eqref{eq:decomposition}.
For this purpose we define
\be
\alpha^k_s = \max_{\bh \in \im\(\bV_{k:0}^s\), \, \| \bh \| = 1} \bh^\T \bP_0 \bh ,
\ee
where $\|. \|$ denotes the Euclidean norm. As a consequence, we have
\be
\(\bV^s_{k:0}\)^\T\bP_0\bV^s_{k:0} \le \alpha^k_s \I{s}.
\ee
From this inequality and from \eqref{eq:decomposition}, we infer
\be
\label{eq:inequality}
\Pi_{\Sv^s}\bM_{k:0}\bP_0 \bM_{k:0}^{\T}\Pi_{\Sv^s} \le  \alpha^k_s \Pi_{\Sv^s}\bU_{k:0}\bSigma^2_{k:0}\bU^\T_{k:0}\Pi_{\Sv^s}.
\ee
Note that, if $\sigma^0_1$ is the largest eigenvalue of $\bP_0$, we have the uniform bound $\alpha^{k}_s \le \sigma^0_1$
for any $k$ and $s$ (see appendix \ref{app:cone}, point 5).  Hence, we can define a finite bound
\be
\alpha_s = \sup_{k \ge 0} \alpha^k_s
\ee
which satisfies for any $k$ and $s$: $\alpha^k_s \le \alpha_s \le  \sigma^0_1$.
Using this uniform bound, in conjunction with \eqref{eq:bound1} and \eqref{eq:inequality}, we obtain
\be
\label{eq:projection}
\Pi_{\Sv^s}\bP_k \Pi_{\Sv^s}  \le \alpha_s \sum_{i\in \Svi^s} \exp\(2\lambda^k_i k\)\bu_i^{k:0} \(\bu_i^{k:0}\)^{\T} .
\ee
Hence, for every unit vector $\bh \in \Sv^s$
\be
\label{eq:ratebound}
\bh^{\T} \bP_k \bh \le  \alpha_s \exp\(2\lambda^k_{n-s+1} k\).
\ee
In particular, if $i \in \Sli$, then $\(\bu^{k:0}_i\)^\T\bP_k\bu^{k:0}_i \rightarrow 0^+$ as $k \rightarrow \infty$.
This defines a {\em weak} form of collapse of $\bP_k$ onto the unstable-neutral subspace $\Ul$.  A {\em strong} form of
collapse is defined by the stable subspace $\Sl$ being in the null space of $\bP_k$.  This can be obtained under the
hypothesis that $\bP_k$ is uniformly bounded, which can in turn be satisfied if the system is sufficiently observed.
Indeed, if $\bP_k$ is uniformly bounded and owing to its positive semi-definiteness, it can be shown that $\left\|
\bP_k\bu^{k:0}_i \right\| \rightarrow 0$ as $k \rightarrow \infty$ (see appendix \ref{app:cone}, point 5).  Hence,
asymptotically, the stable subspace $\Sl$ is in the null space of $\bP_k$. As described in the introduction, this
property is at the core of the class of DA algorithms referred to as AUS \cite[and references therein]{palatella2013a}.

\subsection{Rate of convergence of the eigenvalues}
\label{sec:rate}

In the case of weak -- a fortiori strong -- collapse, the rate of convergence of each of the eigenvalues of $\bP_k$
can be determined from \eqref{eq:ratebound} as follows. Let $\sigma_i^k$ for $i = 1, \dots, n$ denote the
eigenvalues of $\bP_k$ ordered as $\sigma_1^k \ge \sigma_2^k \dots \ge \sigma_n^k$.  Equation~\eqref{eq:ratebound}
guarantees that (appendix \ref{app:cone}, point 6) $\bP_k$ has at least $\s$ of its eigenvalues less than or equal to
$\alpha_s \exp\(2\lambda^k_{n-s+1} k \)$. It follows that
\be
\label{eq:convrate}
\sigma_{i}^k \leq \alpha_i \exp\(2\lambda^k_{i} k \)
\ee
which gives us an upper bound for all eigenvalues of $\bP_k$ and a rate of convergence for the $n-n_0$ smallest ones.

\subsection{Asymptotic rank of the error covariance matrix}
\label{sec:asrank}

A consequence of \eqref{eq:projection} is the upper bound of the asymptotic rank of the error
covariance matrix $\bP_k$. In fact, the asymptotic rank of $\bP_k$ is bounded by the minimum between the rank of $\bP_0$
and $n_0$. This mathematically reads
  \be
  \label{eq:minrank}
  \lim\limits_{k \rightarrow \infty}\rank(\bP_k) \le \min\left\{ \rank(\bP_0), n_0 \right\} .
  \ee

\subsection{Observability and boundedness of the error statistics}
\label{sec:ObsBound}

As mentioned in section \ref{sec:PkP0}, we define the system to be observable if $\mathrm{det}(\bGamma_k) \neq 0$ or,
equivalently, given that $\bM_{k:0}$ is assumed to be non-singular, $\mathrm{det}(\bTheta_k) \neq 0$
\cite{jazwinski1970}. If the system is observable, the inequalities \eqref{eq:bound1} and \eqref{eq:bound2} can
be combined to obtain
\be
\label{eq:bound3}
\bP_k \le \min \left\{ \bM_{k:0}\bP_0 \bM_{k:0}^\T, \bGamma_k^{-1}\right\} ,
\ee
where the accurate definition of the minimum in $\cone{n}$ is given in appendix \ref{app:cone} (point 4).  We note that
if $\bGamma_k$ is bounded by $\bL$ in $\conep{n}$, $\bGamma_k \ge \bL$, we have $\bP_k \le \bGamma_k^{-1} \le \bL^{-1}$
which bounds the error covariances.  The existence of the bound $\bL$ in $\conep{n}$ guarantees the observability of the
system; it forces the precision of the observations to be spread in space and time.  Interestingly, the inequality
\eqref{eq:bound3} reveals that the uncertainty in the state estimate cannot exceed that associated with the most
precise ingredient of the assimilation, the forecast initial conditions or the observations. This is further explored in
the following section.

\section{Asymptotic behavior of $\bP_k$ and its independence from $\bP_0$}
\label{sec:asymptote}

In this section, we study the asymptotic behavior of the forecast error covariance matrix $\bP_k$ when $k \rightarrow
\infty$. In particular, we are interested in the conditions for which the asymptotic sequence of $\bP_k$ becomes
independent of $\bP_0$.  The authors of \cite{pham1998} have provided an appealing and yet heuristic derivation of the
asymptotic limit of $\bP_k$ in the autonomous case, under some observability condition and assuming the absence of a
neutral mode in the dynamics, but also assuming the non-degeneracy of the eigenspectrum of the dynamics.  Here, we are
interested in a rigorous, non-autonomous generalization in the possible presence of neutral modes using a generalized
observability condition. As in the previous section, we assume that the Lyapunov spectrum of the dynamics is
non-degenerate, i.e., there are $n$ distinct Lyapunov exponents. The degenerate case will be discussed at a more
heuristic level at the end of the section.

If $\bC_k$ is a matrix in $\R{n}{n}$ whose columns are the normalized-to-one covariant Lyapunov vectors (CLVs)
of the dynamics at $t_k$ we have the defining relationship
\be
\label{eq:clvdef}
\bM_{k:l}\bC_l = \bC_k\La{k:l},
\ee
where $\La{k:l}$ is a diagonal matrix because of the non-degeneracy of the Lyapunov spectrum. Its diagonal entries are
the exponential of the local Lyapunov exponents between $t_l$ and $t_k$. We will however distinguish between
$\La{k:l}$ and $\La{k:l}^\T$ as if the matrix was not symmetric to ease the discussion on the degeneracy case.  Assuming
that the columns of $\bC_k$ are ordered according to the associated decreasing Lyapunov exponents, we can decompose $\bC_k$
into $\left[ \begin{array}{cc} \bC_{+,k} & \bC_{-,k} \end{array} \right]$ where $\bC_{+,k}$ contains the unstable and
neutral CLVs and $\bC_{-,k}$ contains the stable CLVs.  The transpose of the inverse of $\bC_k$ which, by construction,
forms a dual basis for the CLVs can be decomposed accordingly:
\be
\wbC_k \triangleq  \bC_k^{-\T}  \triangleq \left[ \begin{array}{cc} \wbC_{+,k} & \wbC_{-,k} \end{array} \right] ,
\ee
where $\wbC_{+,k} \in \R{n}{n_0}$ and $\wbC_{-,k} \in \R{n}{(n-n_0)}$.
We decompose $\La{k:l}$ into
\be
\La{k:l} \triangleq
\mat{\Lp{k:l}}{\Lm{k:l}} ,
\ee
where $\Lp{k:l} \in \R{n_0}{n_0}$ and $\Lm{k:l} \in
\R{(n-n_0)}{(n-n_0)}$. Thus, one has
\be
\label{eq:model-split}
\bM_{k:l} = \bC_{+,k}\Lp{k:l}\wbC_{+,l}^\T + \bC_{-,k}\Lm{k:l}\wbC_{-,l}^\T .
\ee
Recall that the FLVs and BLVs are the columns of $\bV_k = \lim_{k \rightarrow \infty} \bV_{k:l}$ and $\bU_k = \lim_{l
  \rightarrow -\infty} \bU_{k:l}$, respectively. Moreover, the FLVs and BLVs associated with the unstable and neutral
directions are the columns of $\bV_{+,k} \in \R{n}{n_0}$ and $\bU_{+,k} \in \R{n}{n_0}$, respectively, which correspond
to the first $n_0$ columns of $\bV_k$ and $\bU_k$, respectively. Finally, the BLVs associated with the stable directions
are the columns of $\bU_{-,k} \in \R{n}{(n-n_0)}$, which correspond to the last $n-n_0$ columns of $\bU_k$.  From
\cite{kuptsov2012}, we have $\bC_k = \bU_k\bT_k = \bV_k\bL_k$, where $\bT_k$ and $\bL_k$ are an invertible upper
triangular matrix and a lower invertible triangular matrix, respectively.  Hence, there is an invertible upper
triangular matrix, $\bT_{+,k} \in \R{n_0}{n_0}$, such that $\bC_{+,k}=\bU_{+,k}\bT_{+,k}$ yielding
$\im(\bC_{+,k})=\im(\bU_{+,k})$.  Moreover, $\wbC_k = \bU_k \bT_k^{-\T} = \bV_k\bL_k^{-\T}$, which implies that there is
an invertible lower triangular matrix and an invertible upper triangular matrix, $\bT^{-\T}_{-,k} \in
\R{(n-n_0)}{(n-n_0)}$ and $\bL^{-\T}_{-,k} \in \R{n_0}{n_0}$, respectively, such that
$\wbC_{-,k}=\bU_{-,k}\bT^{-\T}_{-,k}$ and $\wbC_{+,k}=\bV_{+,k}\bL^{-\T}_{+,k}$. Consequently,
$\im(\wbC_{-,k})=\im(\bU_{-,k})$ and $\im(\wbC_{+,k})=\im(\bV_{+,k})$. These identities will be used in the following.

An asymptotic sequence $\bS_k$ such that $\lim_{k \rightarrow \infty} \(\bP_k - \bS_k\)=\bzero$ will be called an
asymptote for $\bP_k$ in the following.  Our goal is to prove that the two following conditions are sufficient for the
existence of an asymptote for $\bP_k$ which is independent of $\bP_0$. An additional condition may be required if
neutral modes are present in the dynamics.  Recall that $\bP_0$ is possibly degenerate of rank $r_0 \le n$. As in
section \ref{sec:bound}, $\bP_0$ can be factorized into $\bP_0 = \bX_0 \bX_0^\T$, where $\bX_0$ is a matrix in
$\R{n}{r_0}$.

\medskip
\paragraph {\bf Condition 1}
The condition reads
\be
\label{eq:condition1}
\rank \(\wbC_{+,0}^\T\bX_0\) = n_0 .
\ee
The idea is to make the column space of $\bP_0$ large enough so that the unstable and neutral CLVs at $t_0$ have
non-zero projections onto this space.  Since we showed that $\im(\wbC_{+,k})=\im(\bV_{+,k})$, the condition is
equivalent to $\rank \( \bV_{+,0}^\T\bX_0 \) = n_0.$ Consequently, the column space of $\bP_k$ will asymptotically
contain the unstable-neutral subspace.  Note that \eqref{eq:condition1} implies $r_0 \ge n_0$, but $r_0 \ge n_0$
does not imply \eqref{eq:condition1}.

\medskip
\paragraph {\bf Condition 2}
The unstable and neutral directions of the model are uniformly observed, i.e., for $k$ large enough there is
$\varepsilon > 0$ such that
\be
\label{eq:condition2}
\bC_{+,k}^\T \bGamma_k \bC_{+,k} > \varepsilon \I{n_0} .
\ee
The condition is equivalent to $\bU_{+,k}^\T \bGamma_k \bU_{+,k} > \varepsilon \I{n_0}$ with a possibly different $\varepsilon>0$,
since we showed that $\im(\bC_{+,k})=\im(\bU_{+,k})$.

\medskip
We would like to project the degrees of freedom in $\bP_0$ onto the unstable-neutral and stable subspaces. Since
$\bC_{+,0}\wbC^\T_{+,0}+\bC_{-,0}\wbC^\T_{-,0}=\I{n}$, we have
\be
\label{eq:decomp}
\bX_0 = \bC_{+,0}\wbC^\T_{+,0}\bX_0 + \bC_{-,0}\wbC^\T_{-,0}\bX_0 .
\ee
Define $\bZ_+ \triangleq \wbC^\T_{+,0}\bX_0 \in \R{n_0}{r_0}$ which is of rank $n_0$ by Condition 1.  The column spaces of
$\bC_{+,0}\bZ_+$ and of $\bC_{-,0}\wbC^\T_{-,0}\bX_0$ are linearly independent and their sum spans the column space of
$\bX_0$. Hence, they are complementary subspaces in $\im\(\bX_0\)$.
That is why $\bC_{-,0}\wbC^\T_{-,0}\bX_0$ must be of dimension $r_0 - n_0$.  Thus, a QR decomposition of rank $r_0-n_0$ can be used:
$\bC_{-,0}\wbC^\T_{-,0}\bX_0 = \bW_{-,0}\bZ_-$, where $\bW_{-,0} \in \R{n}{(r_0-n_0)}$ is an orthonormal matrix such that
$\bW_{-,0}^\T\bW_{-,0} = \I{r_0-n_0}$ and $\bZ_- \in \R{(r_0-n_0)}{r_0}$ is a full-rank matrix. Note that
$\wbC^\T_{+,0}\bW_{-,0} = \bzero$.  Hence, we have
\be
  \bX_0 = \bC_{+,0}\bZ_+ + \bW_{-,0}\bZ_- 
  = \left[ \begin{array}{cc} \bC_{+,0} & \bW_{-,0} \end{array} \right]\bZ, \quad \mathrm{where} \quad
  \bZ \triangleq \left[ \begin{array}{c} \bZ_+ \\ \bZ_- \end{array} \right] 
\ee
in $\R{r_0}{r_0}$ is of rank $r_0$ since $\rank\(\bX_0\)=r_0$, and hence $\bZ$ is invertible.   
Let us define
\be
\bW_0 \triangleq \left[ \begin{array}{cc} \bC_{+,0}  & \bW_{-,0} \end{array} \right]
\quad \mathrm{and} \quad
\bG \triangleq \bZ^{-\T}\bZ^{-1} ,
\ee
where $\bG \in \conep{r_0}$. Thus $\bP_0 = \bW_0 \bG^{-1} \bW_0^\T$.
This factorization is applied to \eqref{eq:final}:
\begin{align}
  \label{eq:factorization}
  \bP_k &=  \bM_{k:0} \bW_0\bG^{-1}\bW_0^\T\( \I{n} + \bTheta_k  \bW_0\bG^{-1}\bW_0^\T\)^{-1} \bM_{k:0}^\T \nn
  &= \bM_{k:0} \bW_0\bG^{-1}\( \I{r_0} + \bW_0^\T\bTheta_k  \bW_0\bG^{-1}\)^{-1} \bW_0^\T\bM_{k:0}^\T \nn
  &= \bM_{k:0} \bW_0\(\bG + \bW_0^\T\bTheta_k  \bW_0\)^{-1} \bW_0^\T\bM_{k:0}^\T ,
\end{align}
where appendix \ref{app:sml} has been employed.

We now focus on the projection of $\bP_k$ onto the forward unstable-neutral subspace, which reads, from
(\ref{eq:thetadef},\ref{eq:model-split},\ref{eq:factorization}),
\begin{align}
  \label{eq:block1}
  \wbC^\T_{+,k}\bP_k\wbC_{+,k}  & =  \Lp{k:0}\wbC^\T_{+,0}\bW_0 \(\bG + \bW_0^\T\bTheta_k  \bW_0\)^{-1}\bW_0^\T \wbC_{+,0} \Lp{k:0}^\T\nn
  & =  \Lp{k:0} \left[ \(\bG + \bW_0^\T\bTheta_k  \bW_0\)^{-1} \right]_{\! _{++}} \Lp{k:0}^\T\nn
  & = \Lp{k:0}
  \left[ \( \begin{array}{cc}
      \Gpp + \Tpp &
      \Gpm + \Tpm \\
      \Gmp + \Tmp &
      \Gmm + \Tmm
    \end{array}\)^{-1}\right]_{\! _{++}}
  \Lp{k:0}^\T \nn
 &  = \left[ \( \begin{array}{cc}
      \Lp{k:0}^{-\T}\Gpp\Lp{k:0}^{-1} + \Gapp &
       \Lp{k:0}^{-\T}\Gpm + \Lp{k:0}^{-\T}\Tpm  \\
      \Gmp\Lp{k:0}^{-1} + \Tmp\Lp{k:0}^{-1}&
      \Gmm + \Tmm
    \end{array}\)^{-1}\right]_{\! _{++}}
\end{align}
where $\left[\cdot \right]_{\! _{++}}$ is the matrix block corresponding to the unstable-neutral modes, and
\begin{align}
  &\Gapp \triangleq \bC_{+,k}^\T\bGamma_k\bC_{+,k} = \Lp{k:0}^{-\T}\Tpp\Lp{k:0}^{-1} ,
  \quad \Tpp \triangleq \bC_{+,0}^\T\bTheta_k\bC_{+,0} ,   \\
  &\Tmp \triangleq \bW_{-,0}^\T\bTheta_k\bC_{+,0} ,
  \qquad  \Tpm \triangleq \bC_{+,0}^\T\bTheta_k\bW_{-,0} , \qquad
  \Tmm \triangleq \bW_{-,0}^\T\bTheta_k\bW_{-,0} .
\end{align}

\subsection{Asymptote in the absence of neutral modes}
\label{sec:asymptote-nn}

We first assume the absence of any neutral mode in the dynamics.  In \eqref{eq:block1}, the term
$\Lp{k:0}^{-\T}\Gpp\Lp{k:0}^{-1}$ asymptotically vanishes since in the absence of neutral modes,
\be
\label{eq:aslim}
\lim_{k \rightarrow  \infty} \left\| \Lp{k:l}^{-1} \right\| = \bzero
\ee
for any matrix norm $\left\| \cdot \right\|$ and any $l$. The behavior of $\bE_{+,k} \triangleq
\Gmp\Lp{k:0}^{-\T}+\Tmp\Lp{k:0}^{-1}$ and $\bE_{+,k}^\T$ as the unstable/stable off-diagonal blocks in \eqref{eq:block1}
remain to be studied. To that end, we choose a submultiplicative matrix norm, and obtain 
\begin{align}
  \left\| \bE_{+,k}  \right\| \le & \left\| \Gmp \right\|\left\| \Lp{k:0}^{-\T} \right\| \nn
  & + \sum_{l=0}^{k-1} \left\| \bW_{-,0}^\T\right\|
\left\|\wbC_{-,0}\right\| \left\| \Lm{l:0}^\T \right\| \left\|\bC_{-,l}^\T \right\|
\left\|\bOmega_l \right\| \left\|\bC_{+,l} \right\| \left\|\Lp{k:l}^{-1} \right\|  .
\end{align}
Since $\bOmega_l$ has been assumed uniformly bounded from above, and given that matrices with unitary columns are uniformly
bounded (by $\sqrt{n}$ in Frobenius norm), there is a constant $c_+$ such that for all $k$
\be
\left\| \bE_{+,k} \right\| \le  \left\| \Gmp \right\|\left\| \Lp{k:0}^{-\T} \right\|
+ c_+ \sum_{l=0}^{k-1}  \left\| \Lm{l:0} \right\| \left\|\Lp{k:l}^{-1} \right\|.
\ee
Here the matrix norm can be thought of as the spectral norm up to a multiplicative constant.  Because the sum is dominated
by the convergent series $\sum_{l=0}^ \infty \left\| \Lm{l:0} \right\|$ and because of \eqref{eq:aslim} for each
$l$, the majorant of $\bE_{+,k}$ asymptotically vanishes so that $\lim_{k \rightarrow \infty} \left\| \bE_{+,k} \right\|
= \bzero$.

As a consequence, the off-diagonal terms of $\bG + \bW_0^\T\bTheta_k \bW_0$ in \eqref{eq:block1} asymptotically
vanish.  Moreover, the diagonal blocks are uniformly bounded from below by $\varepsilon \I{n_0}$ by Condition 2 for the
top-left block and $\Gmm \in \conep{n-n_0}$ for the bottom-right block.  Consequently, the inverse of $\bG +
\bW_0^\T\bTheta_k \bW_0$ is asymptotically given by the inverse of the diagonal blocks, up to a vanishing sequence of
matrices, and for the unstable block one obtains
\be
\lim_{k \rightarrow \infty} \left\{ \wbC^\T_{+,k}\bP_k\wbC_{+,k}- \( \bC_{+,k}^\T\bGamma_k\bC_{+,k} \)^{-1}\right\} = \bzero .
\ee
This implies by Condition 2, that $\wbC^\T_{+,k}\bP_k\wbC_{+,k}$ is asymptotically uniformly bounded.
Moreover, we have proven in section \ref{sec:PkCONS} that $ \lim_{k \rightarrow \infty} \bU^\T_{-,k}\bP_k\bU_{-,k} =
\bzero$, which by $\im(\wbC_{-,k})=\im(\bU_{-,k})$, shows that $ \lim_{k \rightarrow \infty} \wbC^\T_{-,k}\bP_k\wbC_{-,k} =
\bzero$\footnote{Alternatively, this can be recovered from \eqref{eq:factorization} including the rate of convergence.}.
As a consequence of appendix \ref{app:cone}, point 5, and using
the fact that $\bC_k^{-1}\bP_k\bC_k^{-\T}$ is in $\cone{n}$, its 
off-diagonal blocks $\wbC^\T_{+,k}\bP_k\wbC_{-,k}$ and $\wbC^\T_{-,k}\bP_k\wbC_{+,k}$ also asymptotically vanish. We conclude
(using the uniform boundedness of the $\bC_k$)
\be
\label{eq:aseq}
\lim_{k \rightarrow \infty} \left\{ \bP_k - \bS_k \right\} = \bzero, \quad \mathrm{where}
\quad \bS_k = \bC_{+,k}\( \bC_{+,k}^\T\bGamma_k\bC_{+,k} \)^{-1} \bC_{+,k}^\T .
\ee
Importantly, the asymptote $\bS_k$ does not depend on $\bP_0$.  This generalizes the result in
\cite{pham1998} to the non-autonomous case. Moreover, we will see later that this also generalizes their result to the
case where the Lyapunov spectrum is degenerate.

The case where neutral modes are present is much more involved, and yet physically very important
\cite{vannitsem2016}. Indeed, there is a wide range of possible asymptotic behaviors for neutral CLVs, which could, for
instance, grow, or decay, at a subexponential rate.  They could intermittently behave as unstable or stable modes.  To
go further in the case where neutral modes are present, and only then, the neutral modes and their observability need to
be characterized more precisely.  Thus, an additional condition that complements Condition 2 needs to be introduced.
Since several conditions are possible, we focus on one of them and discuss two others.

\subsection{Asymptote in the presence of neutral modes}

It is convenient to further split the unstable from the neutral CLVs as
\be \bC_k \triangleq \left[ \begin{array}{ccc} \bC_{++,k} & \bC_{+0,k} & \bC_{-,k} \end{array} \right] .
\ee
Accordingly, we refine the decomposition
\be
\La{k:l} \triangleq
\left[ \begin{array}{ccc}
    \Lpp{k:l} & \bzero   & \bzero \\
    \bzero   & \Lpn{k:l} & \bzero \\
    \bzero   & \bzero   & \Lm{k:l}
  \end{array} \right] ,
\ee
where $\Lpp{k:l} \in \R{n_{0+}}{n_{0+}}$ and $\Lpn{k:l} \in \R{n_{00}}{n_{00}}$ are diagonal with $n_{0+}+n_{00}=n_0$.
Furthermore, $\Lpn{k:l}$ is factorized into growing and decaying contributions.
First, $\Lpng{k:0} \in \R{n_{00}}{n_{00}}$ is a diagonal matrix with entries: for $i=1,\ldots,n_{00}$,
\be
   [\Lpng{k:l}]_{ii} = [\Lpn{k:l}]_{ii} \prod_{q=l+1}^k \max\(\left|[\Lpn{q:q-1}]_{ii}\right|^{-1},1\) .
\ee
This definition performs the factorization while transferring the transitivity property of $\La{k:l}$, i.e.,
$\La{k:0}=\La{k:l}\La{l:0}$ for $k \ge l \ge 0$, to $\Lpng{k:l}$, i.e., $\Lpng{k:0}=\Lpng{k:l}\Lpng{l:0}$ for $k \ge l
\ge 0$.  Second, we also define $\Lpnl{k:0} \triangleq \Lpn{k:0}\Lpng{k:0}^{-1}$ and
\be
\Lpg{k:0} \triangleq \mat{\Lpp{k:0}}{\Lpng{k:0}} ,
\quad 
\Lpl{k:0} \triangleq \mat{\I{n_{0+}}}{\Lpnl{k:0}}.
\ee
Let us suggest a possible condition that supplements Conditions 1 and 2, in order to obtain an asymptote $\bS_k$ for
$\bP_k$ which, in the presence of neutral modes, does not depend on $\bP_0$.

\medskip
\paragraph {\bf Condition 3} We define 
\be
\label{eq:phidef}
\bPhi_{k} =  \Lpng{k:0}^{-\T}\bC^\T_{+0,0}\bTheta_k\bC_{+0,0} \Lpng{k:0}^{-1}
\ee
and require that for any $\bv \in {\mathbb R}^{n_{00}}$, and using the Euclidean norm $\left\| \cdot\right\|$
\be
\label{eq:condition3}
\liminf\limits_{k \rightarrow \infty} \left\| \bPhi_{k} \bv \right\| = +\infty .
\ee
Let us give the example of a realistic setup that leads to Condition 3.  Assume that the dynamics is autonomous. Hence,
there is at least one neutral CLV $\bv$ for which the local Lyapunov exponents are $0$.  Suppose that there is a
sequence $\omega_l$ of non-negative numbers such that $\bOmega_l \ge \omega_l \bv \bv^\T$, and assume $\sum_{l=0}^\infty
\omega_l = \infty$, then it is not difficult to show through (\ref{eq:thetadef},\ref{eq:phidef}) that Condition 3
is satisfied. In particular, if the observation process is time-invariant such that $\omega \triangleq \omega_l > 0$,
then $\left\| \bPhi_{k} \bv \right\| = \omega k\left\|\bv \right\|$ which indeed diverges with $k$.

Analysis of \eqref{eq:block1} in the presence of neutral modes is not straightforward, and instead we use a Schur
complement for the inverse of the top-left block,
\begin{align}
  \label{eq:block2}
  \left[ \(\bG + \bW_0^\T\bTheta_k  \bW_0\)^{-1} \right]^{-1}_{\! _{++}}
  = & \Gpp+\Tpp -  \(\Gpm+\Tpm\)\times \nn
  & \(\Gmm+\Tmm\)^{-1}\(\Gmp+\Tmp\) .
\end{align}
Separating the growing and decaying trends of the neutral modes, this yields
\begin{align}
\label{eq:block3}
\Lpl{k:0}^\T\( \wbC^\T_{+,k}\bP_k\wbC_{+,k}\)^{-1}\Lpl{k:0}
=& \Lpg{k:0}^{-\T} \(\Gpp + \Tpp\) \Lpg{k:0}^{-1} \nn
& - \bE_k^\T \(\Gmm+\Tmm\)^{-1}\bE_k,
\end{align}
where $\bE_k \triangleq \Gmp\Lpg{k:0}^{-1} + \Tmp\Lpg{k:0}^{-1}$. $\bE_k$ can be split into $\bE_k =
\left[ \begin{array}{cc }\bE_{+,k} & \bE_{0,k} \end{array} \right]$, where $\bE_{+,k} \triangleq
\bW_{-,0}^\T\bG\bC_{++,0}\Lpp{k:0}^{-1} + \bW_{-,0}^\T\bTheta_k\bC_{++,0}\Lpp{k:0}^{-1}$ is the stable/unstable matrix
block and $\bE_{0,k} \triangleq \bW_{-,0}^\T\bG\bC_{+0,0}\Lpng{k:0}^{-1} +
\bW_{-,0}^\T\bTheta_k\bC_{+0,0}\Lpng{k:0}^{-1}$ is the stable/neutral matrix block.  This definition of $\bE_{+,k}$ is
consistent with that of $\bE_{+,k}$ in section \ref{sec:asymptote-nn}.  For the same reasons as in section
\ref{sec:asymptote-nn}, we have $\lim_{k \rightarrow \infty} \left\| \bE_{+,k} \right\| = \bzero$.  Similarly, there is
a constant $c_0$ such that for all $k$
\be
\label{eq:bound4}
\left\| \bE_{0,k} \right\| \le \left\| \Gmp \right\|\left\| \Lpng{k:0}^{-1} \right\| + c_0 \sum_{l=0}^{k-1}
\left\| \Lm{l:0} \right\| \left\|\Lpnl{l:0} \right\|\left\|\Lpng{k:l}^{-1} \right\|,
\ee 
which, however, only ensures that $\left\|\bE_{0,k}\right\|$ is uniformly bounded from above. Note that, in deriving
\eqref{eq:bound4}, we used
$\Lpn{l:0}\Lpng{k:0}^{-1}=\Lpn{l:0}\Lpng{l:0}^{-1}\Lpng{k:l}^{-1}=\Lpnl{l:0}\Lpng{k:l}^{-1}$.  Hence, the only case
where the last term of \eqref{eq:block3}, as well as $\Lpg{k:0}^{-\T} \Gpp \Lpg{k:0}^{-1}$, do not asymptotically
vanish but are uniformly bounded is for entries related to the neutral modes alone.  Yet, in this case, the neutral
sub-block of $\Lpg{k:0}^{-\T} \Tpp \Lpg{k:0}^{-1}$, which is $\bPhi_k$, asymptotically dominates the bounded correction
by Condition 3, so that, even in the presence of neutral modes, the correction term is negligible.  More precisely,
\eqref{eq:block3} has the asymptote
\be
\left[ \begin{array}{cc}
    \bC_{++,0}^\T\bGamma_k\bC_{++,0} & \Lpp{k:0}^{-\T}\bC_{++,0}^\T\Tpp\bC_{+0,0}\Lpng{k:0}^{-1} \\
    \Lpng{k:0}^{-\T}\bC_{+0,0}^\T\Tpp\bC_{++,0}\Lpp{k:0}^{-1} & \bPhi_k + \bB_k
\end{array} \right] ,
\ee
where $\bB_k$ is a bounded sequence. Its inverse, i.e., $\Lpl{k:0}^{-1} \wbC^\T_{+,k}\bP_k\wbC_{+,k}\Lpl{k:0}^{-\T}$,
is well-defined and obtained by the formula of the inverse of a matrix with $2\times 2$ sub-blocks, which yields the asymptote
\be
\left[ \begin{array}{cc}
    \( \bC_{++,0}^\T\bGamma_k\bC_{++,0} \)^{-1} & \bzero \\
    \bzero & \bPhi_k^{-1}
\end{array} \right] .
\ee
We can conclude similarly to the end of section \ref{sec:asymptote-nn} that \eqref{eq:aseq} is still valid in the
presence of neutral modes if Condition 3 is satisfied.  In this case, the asymptotic sequence can further be simplified
into
\be
\label{eq:aseq2}
\bS_k = \bC_{++,0} \( \bC_{++,0}^\T\bGamma_k\bC_{++,0} \)^{-1} \bC_{++,0}^\T .
\ee
Note that we expect the rate of convergence to the neutral subspace to be possibly quite different from the exponential
convergence to the stable subspace. For instance, considering again the example where $\left\| \bPhi_{k} \bv_k \right\|
= \omega k\left\|\bv_k \right\|$, $\wbC_{+0,k}^\T\bP_k\wbC_{+0,k}$ asymptotically behaves like $(\omega
k)^{-1}\I{n_{00}}$.

An alternative to Condition 3, but which does not exhaust all the possible behaviors of neutral dynamics, is to
assume that the sequence $\Lpng{k:0}^{-1}$ has a limit and
\be
\lim_{k \rightarrow \infty} \left\| \Lpng{k:0}^{-1} \right\| = \bzero .
\ee
If satisfied, $\bE_{0,k}$ and $\Lpg{k:0}^\T\Gpp\Lpg{k:0}^{-1}$ in \eqref{eq:block3} both asymptotically
vanish. Therefore, the convergence is similar to that of section \ref{sec:asymptote-nn}: the neutral modes are
effectively unstable modes. Equation~\eqref{eq:aseq} remains valid but the neutral modes may have a non-negligible
contribution to $\bS_k$ depending on the asymptotic behavior of $\bGamma_{00,k} =
\bC_{+0,k}^\T\bGamma_{k}\bC_{+,0;k}$.

If the sequence $\left\|\Lpng{k:0}^{-1}\right\|$ does not have a limit or if its limit is distinct from $0$, then the
sequence has an adherence point in $]0,1]$.Then, the sequence $\left\| \bE_{0,k} \right\|$ does not necessarily
    asymptotically vanish.  However, it can be seen that if the neutral modes are observed with an $\bOmega_l$ bounded
    from below, $\bPhi_k$ necessarily diverges, which leads back to the asymptote \eqref{eq:aseq2}.

Note that we have discussed criteria applying to all neutral modes. However, Condition 3 or its alternatives could be
individualized to each neutral mode.

\subsection{Degeneracy of the Lyapunov spectrum}

At a more heuristic level, we discuss the case where there are multiplicities greater than one in the Lyapunov spectrum.
In this case, \eqref{eq:clvdef} remains valid, up to some qualifications. While the Oseledec subspaces themselves
are covariant, they may not be entirely decomposable into one dimensional covariant subspaces -- one should immediately
consider the analogy with generalized eigenvectors for Jordan blocks. Indeed, when there is a degenerate Lyapunov
spectrum there may only be a single covariant vector per exponent. If we redefine $\bC_k$ to be a matrix with columns
composed of an ordered basis for each covariant Oseledec subspace then $\La{k:l}$ is not necessarily a diagonal matrix
but rather upper triangular over $\mathbb C$.  In that case, the transpose operators in the above derivations should be
understood as conjugate and transpose operators.  This consideration only matters in the derivation where norms of
$\La{k:l}$ are to be computed when studying the convergence of sequences and sums. In this case, the spectral norm is
not necessarily related to the eigenspectrum of $\La{k:l}$ which is not necessarily self-adjoint (something that has
been overlooked in the derivation in \cite{pham1998} in the autonomous case).  However, we can heuristically expect a
polynomial correction as a function of $k-l$ to the exponential growth or decay of $\La{k:l}$ because of the impact of
the Jordan blocks.  By the Oseledec theorem, ones knows that $\lim_{k-l \rightarrow \infty}
\(\La{k:l}^\T\La{k:l}\)^{1/2(k-l)} = \exp(\bD)$ where $\bD$ is the diagonal matrix of the Lyapunov exponents.  Hence,
$\La{k:l}^\T\La{k:l}$ asymptotically behaves like $\exp((k-l)\bD)$ up to a subexponential correction.  To compute the
matrix norm of $\La{k:l}$, one can equivalently use the spectral norm and compute the square root of the spectral radius
$\rho(\La{k:l}^\T\La{k:l})$ which is asymptotically equivalent to $\rho(\exp((k-l)\bD))$ up to a subexponential
correction which could be absorbed into the exponential trend if need be.  This shows that the exponential trends in the
above derivations are unchanged.  As a consequence, we believe that the main results of this section are likely to be
unaltered by the Lyapunov spectrum degeneracy. Only the neutral modes, that are not subject to exponential growth or
decay, could be impacted.  This will be studied and reported in a separate paper.

\subsection{Role of the neutral modes}

It is finally worth mentioning the particular role played by the neutral modes.  In section \ref{sec:PkCONS}, the
exponential convergence to $\bzero$ of $\bP_k$ for all stable directions was proven.  Provided the three conditions
(\ref{eq:condition1},\ref{eq:condition2},\ref{eq:condition3}) are met, the above discussion points to an
exponential convergence of $\bP_k$ onto $\bS_k$ for all unstable directions. Nevertheless, the discussion also suggests
a slower convergence of $\bP_k$ to $\bzero$ for neutral directions.  The critical importance of the neutral modes was
originally observed by \cite{trevisan2010} in numerical experiments with assimilation in the unstable-neutral subspace
in the context of variational DA with nonlinear dynamics. They numerically showed that it was necessary
to include the neutral direction within the subspace where the assimilation was performed in order to efficiently
control the error growth. The analysis carried out in the present section further corroborates their findings.

Moving away from the linear hypothesis towards nonlinear dynamics and in connection with this slow convergence of the
neutral modes, it was recently argued \cite{bocquet2015} that the region of the Lyapunov spectrum around the neutral
modes is critical in the convergence of the EnKF. The mis-estimation of the uncertainty in this region of the spectrum
was shown to be the reason why the ad hoc technique known as inflation meant to stabilize the filter is very often
required.

\section{Numerical results}
\label{sec:numeric}

We present here numerical results on the asymptotic properties of the analysis error covariance $\bP^\mathrm{a}_k$ that
corroborate and illustrate the theoretical findings. The convergence results obtained for $\bP_k$ can easily be
transferred to $\bP^\mathrm{a}_k$ by \eqref{eq:kfforecast}, or by applying $\bP^{\rm a}_k \leq \bP_k$ which is
readily obtained using the matrix shift lemma to \eqref{eq:kfupdate} as in \eqref{eq:recurrence2}.

Three different experimental setups are considered, with different choices of the dynamical and observational models in
(\ref{eq:dynmodel},\ref{eq:obsmodel}).  In all cases the perfect model hypothesis is employed, $\bQ_k={\bf 0}$.

\medskip
\paragraph {\bf Exp1: Autonomous system}
The state- and observation-space dimensions are $n=30$ and $d=10$, respectively. The time-invariant matrices $\bM_k
\triangleq \bM \in \mathbb{R}^{n \times n}$, $\bH_k \triangleq \bH\in \mathbb{R}^{d \times n}$, and $\bR_k \triangleq
\bR \in \mathbb{R}^{d \times d}$ are chosen randomly, i.e., with entries which are independently and identically
distributed (iid) standard normal random variables.

\medskip
\paragraph {\bf Exp2: Random non-autonomous system}
The state- and observation-space dimensions are $n=30$ and $d=10$ respectively.  The time-varying, invertible, propagators
$\bM_k \in \mathbb{R}^{n \times n}$, the observation error covariance matrices $\bR_k \in \mathbb{R}^{d \times d}$
and the observation matrices $\bH_k \in \mathbb{R}^{d \times n}$ are all randomly generated, i.e., the entries
of these matrices are iid standard normal random variables.

\medskip
\paragraph {\bf Exp3: Non-autonomous system obtained by linearization around a trajectory of the Lorenz-95 model}
  The entries of the observation error covariance, $\bR_k$ and $\bH_k$ are generated as in Exp2 but with the state- and
  observation-space dimensions being $n=40$ and $d=15$, respectively. The propagators $\bM_k$ are taken to be the
  linearization around a trajectory on the attractor of the $n=40$-dimensional Lorenz-95 model \cite{lorenz1998}, which
  is very commonly used in DA literature; see e.g., \cite[and references therein]{carrassi2009}.
The equations read
\begin{equation}
\frac{{\rm d}x_j}{{\rm d}t} = x_{j-1}\left(x_{j+1} -x_{j-2}\right)-x_j+F, \qquad j=1,...,n  ,
\end{equation}
with periodic boundary conditions, $x_0=x_n$, $x_{-1}=x_{n-1}$, and $x_{n+1}=x_{1}$. The standard value of the forcing,
$F=8$, is used in the following experiments. The observation interval is $\Delta t = 0.1$.

In another numerical experiment, we used a simpler observational network, by choosing $\bH_k = [1, 0, \dots, 0]$,
corresponding to observation of only the first component of the state vector. The numerical results for Exp3 (with
randomly chosen elements of $\bH_k$ of dimension $15 \times 40$) and for this much simpler observational network were
qualitatively the same and thus the latter have not been presented here.

It must be emphasized that this case (Exp3) does not coincide with the nonlinear filtering problem of the Lorenz-95
model but it makes use of the linearization of the model to build up the propagator which is then used as a linear model
in \eqref{eq:dynmodel}.

Note, furthermore, that in an extended Kalman filter (EKF, \cite{jazwinski1970}) applied to a nonlinear system such as
the Lorenz-95, the only place where the state estimate enters the computation of the covariance matrices is in the
linearization of the model dynamics in which one needs to estimate the Jacobian of the model dynamics evaluated on the
system's state. Therefore, for the Lorenz-95 model, the analysis and forecast covariances of the EKF will show asymptotic
behavior similar to what is presented below. While this behavior was already observed and exploited in a reduced-order
formulation of the EKF based on the unstable subspace \cite{trevisan2011}, it does not give many hints about the
asymptotic behavior of a fully nonlinear filter.

Each of the three experimental setups is representative of a class of systems. Numerical results (not shown) for other
choices of the system and observational dimension as well as for other realizations of the random matrices $\bM_k,
\bH_k, \bR_k$ were found to be qualitatively equivalent to the results reported below.

For the three numerical experiments described above, it is very difficult to check the observability condition
$\mathrm{det}(\bGamma_k) \neq 0$ because the matrices $\bGamma_k$ soon become very ill-conditioned. But we expect that
the system would be observable with probability 1 for Exp1 and Exp2, while in Exp3 (the case of Lorenz-95 model),
we expect the system to be observable even with a single variable being observed, since each variable is coupled to
those around it.  Note that we are unable to numerically verify the above statements.

The QR method \cite{parker1989,legras1996} is adopted to numerically compute the Lyapunov vectors and
exponents. Starting from a random positive semi-definite $\bP^{\rm a}_0$, the sequence $(\bP_k, \bP^{\rm a}_k)$ for $k >
0$ of forecast and analysis error covariance matrices was generated based on the KF (\ref{eq:kfupdate},\ref{eq:kfforecast}).

Recall that $n_0$ stands for the number of non-negative Lyapunov exponents: in most cases, $n_0$ will correspond to the
number of positive plus one zero exponent.  Numerically, this zero Lyapunov exponent will not be exactly zero but it
will fluctuate around it.  Also recall that $r_0$ is the rank of the initial covariance matrices $\bP_0$, or $\bP^{\rm a}_0$.

\subsection{Rate of convergence of the eigenvalues}
\label{sec:rate_numer}

The following numerical experiments show the relation between the rates of convergence of eigenvalues, $\sigma_i^k$, of
the error covariance matrix $\bP_k$, and the Lyapunov exponents of the dynamical system of
(\ref{eq:dynmodel},\ref{eq:obsmodel}).  The eigenvalues are ordered so that $\sigma_1^k \ge \sigma_2^k \dots \ge
\sigma_n^k $.

When $r_0 < n_0$, the rank of $\bP_k$ as $k \rightarrow \infty$ generically remains $r_0$ for almost all initial
conditions with no eigenvalues decaying to zero.  Thus we consider here the relevant situation: $r_0 \ge n_0$. From
section \ref{sec:RANKP}, we know that in this case, $\sigma_{r_0+1}^k = \dots = \sigma_{n}^k = 0$.  Moreover, from
section \ref{sec:asrank}, \eqref{eq:minrank}, we know that $\sigma_{1}^k, \dots, \sigma_{n_0}^k$ will remain non-zero
even in the limit $k \to \infty$ -- except maybe for the neutral directions as discussed in section \ref{sec:asymptote} --
whereas $\sigma_{n_0+1}^k, \dots, \sigma_{r_0}^k$ will decay to zero.  Recall from
section \ref{sec:collapse} that $\exp(\lambda_i^k k )$ is a singular value of $\bM_{k:0}$ so that $\lambda_i^k$
approaches the Lyapunov exponent $\lambda_i$ as $k \to \infty$; the Lyapunov exponents are ordered so that the first
$n_0$ are non-negative, $\lambda_1 > \lambda_2 > \dots > \lambda_{n_0} \ge 0$, whereas the rest are negative with
decreasing value, $0 > \lambda_{n_0+1} > \lambda_{n_0+2} > \dots > \lambda_n$.

The results of section \ref{sec:rate} and the inequality \eqref{eq:convrate} can be used to derive 
the rate of convergence of the smallest $r_0-n_0$ eigenvalues $\sigma_{i}^k$ with $i = n_0 + 1, \dots, r_0$. 
Using the largest eigenvalue at the initial time, $\sigma_1^0$, as for $\alpha_i$ in \eqref{eq:convrate} we have
\be
\label{eq:expdecayrates} 
\sigma_{i}^k \le \sigma_1^0 \exp\left(2 \lambda_i^k k  \right) 
\ee
which implies that, asymptotically,
\be
\label{eq:decayrates} 
\ln(\sigma_{i}^k) \leq  \ln(\sigma_1^0) + 2 \lambda^k_i k \underset{k \rightarrow \infty}{\sim}
\ln(\sigma_1^0) + 2 \lambda_i k .
\ee
The equivalence in \eqref{eq:decayrates} is valid in the limit $k \to \infty$ since
$\lambda_i^k \to \lambda_i$ as $k \to \infty$.  Therefore for $i = n_0 + 1, \dots, r_0$, the eigenvalues $\sigma_{i}^k$
of $\bP_k$ decay to zero exponentially fast with the exponential decay rate asymptotically being at least twice the
Lyapunov exponent $\lambda_i$. Note that, as mentioned above, $\bP^{\rm a}_k \leq \bP_k$, so the
aforementioned decay rate is also valid for the eigenvalues of the analysis error covariance matrix, $\bP^{\rm a}_k$.

\begin{figure}[htbp]
  \begin{tabular}{cc}
    \includegraphics[width=0.46 \textwidth]{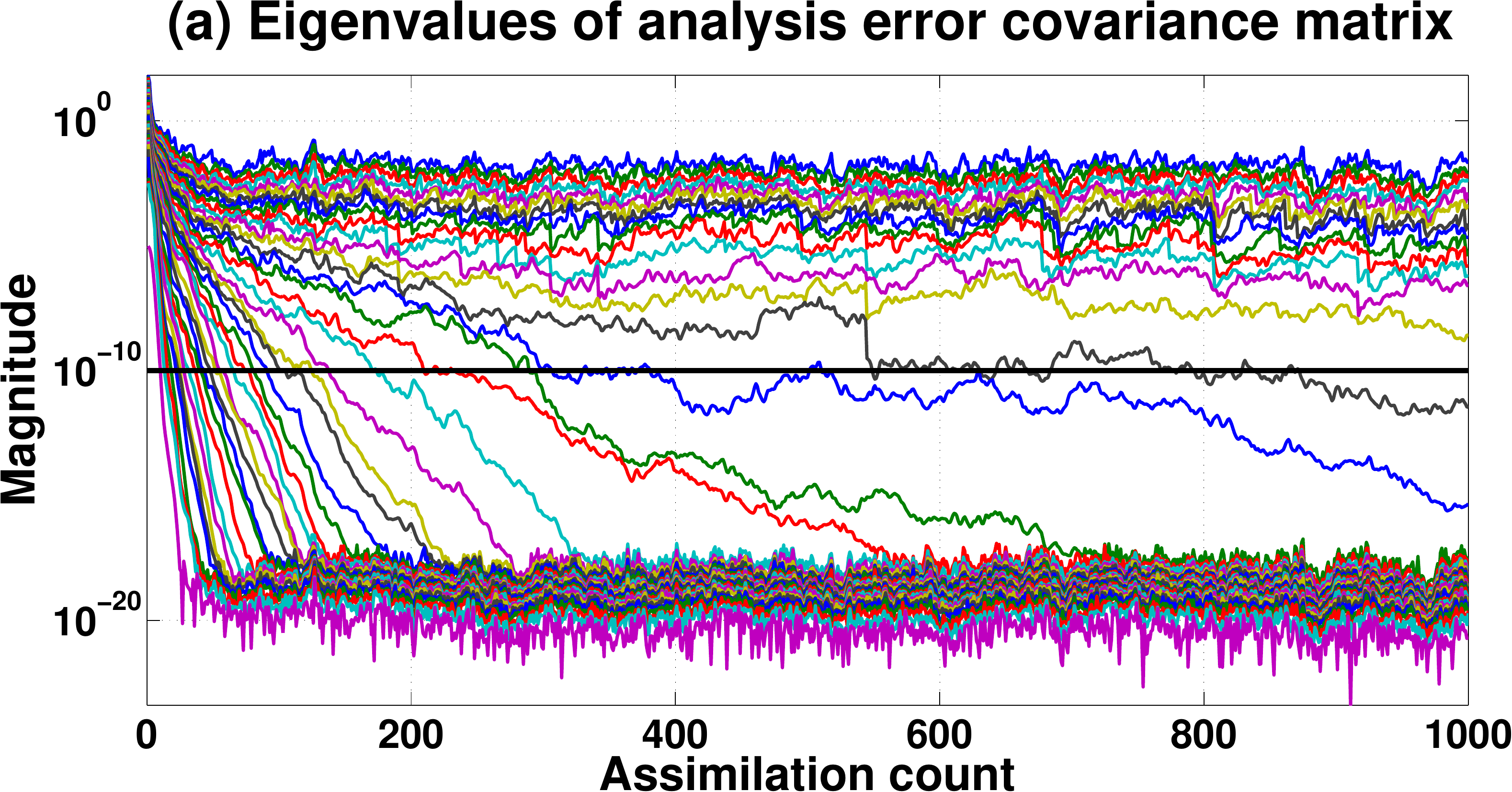} &
    \includegraphics[width=0.46 \textwidth]{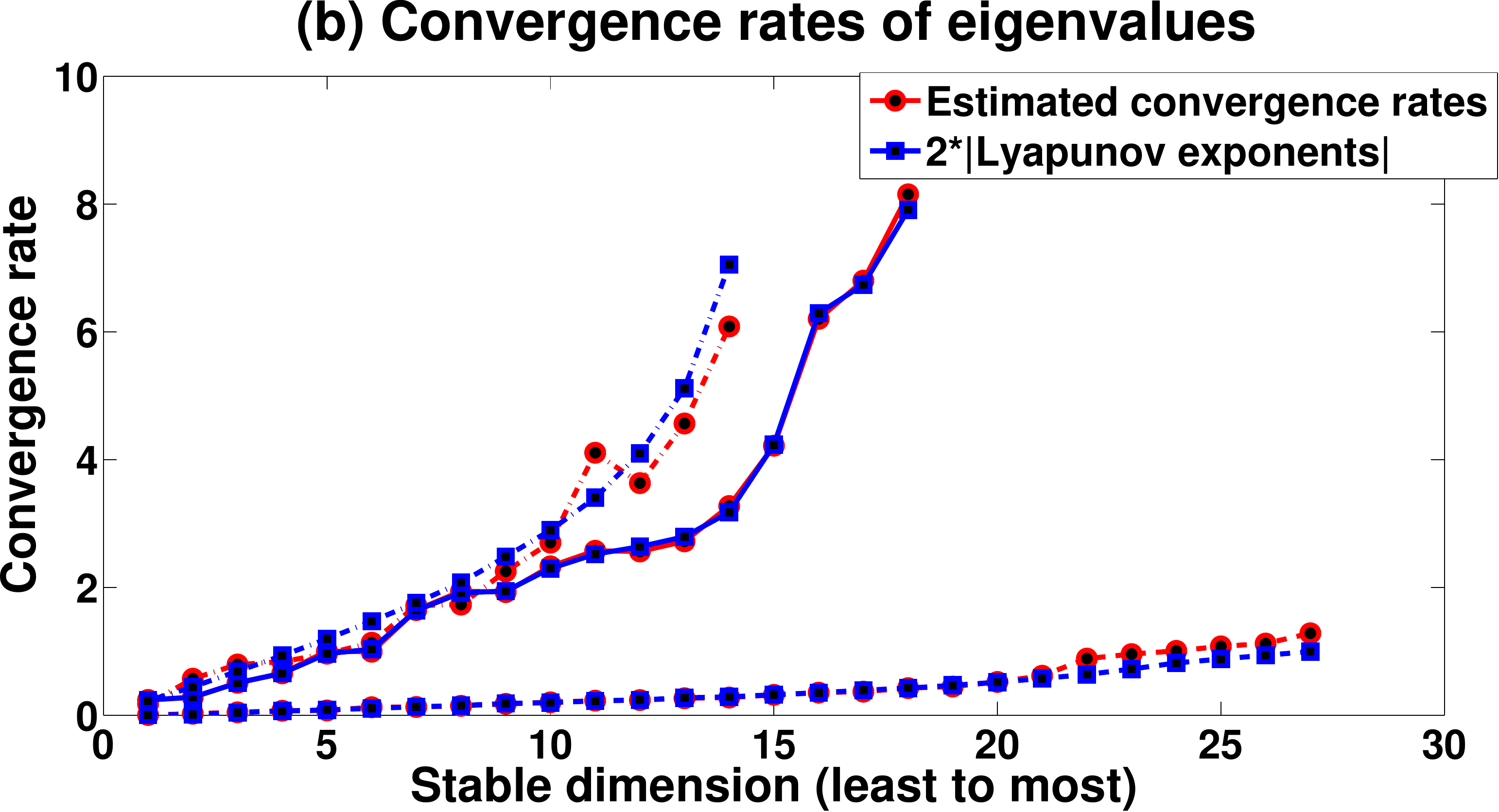}
  \end{tabular}
    \caption{\label{fig:5decay} Panel (a) shows the eigenvalues of $\bP_k^\mathrm{a}$ in Exp3 and the decay of part of
      its spectrum. Panel (b) shows the comparison between the decay rate of the eigenvalues of the analysis covariance
      matrix (red lines) with twice the absolute value of the negative Lyapunov exponents (blue lines), for the
      autonomous system (solid line, $n = 30, n_0=16$), and for two examples of non-autonomous systems with random
      propagators (dash-dot line, $n = 30, n_0=16$) and with propagators derived from the Lorenz-95 system (dashed line,
      $n = 40, n_0=14$), for full rank $\bP_0$.}

\end{figure}

Figure \ref{fig:5decay}(a) illustrates the decay of some of the eigenvalues of the analysis error covariance matrix in
Exp3. Similar graphs have been obtained for Exp1 and Exp2 (not shown).  Figure \ref{fig:5decay}(b) shows the slopes of
the best fit lines for the semi-log plot of $\sigma_{i}^k$ versus $k$ for the full rank $\bP_0$ (red
dots), for all three experiments described above.  The blue dots show the values of twice the absolute value of the
corresponding negative Lyapunov exponents.  We see that, indeed, the inequality \eqref{eq:decayrates} is saturated.

\subsection{Existence of asymptotic sequences of low-rank covariance matrices}
\label{sec:nonaut}

The next set of numerical results corroborate the results in section \ref{sec:collapse} about the projections of the
covariance matrices onto the stable subspace vanishing and the results in section \ref{sec:asymptote} about their
asymptotic behavior.

\begin{figure*}[htbp]
  \begin{tabular}{cc}
    \includegraphics[width=0.46\textwidth]{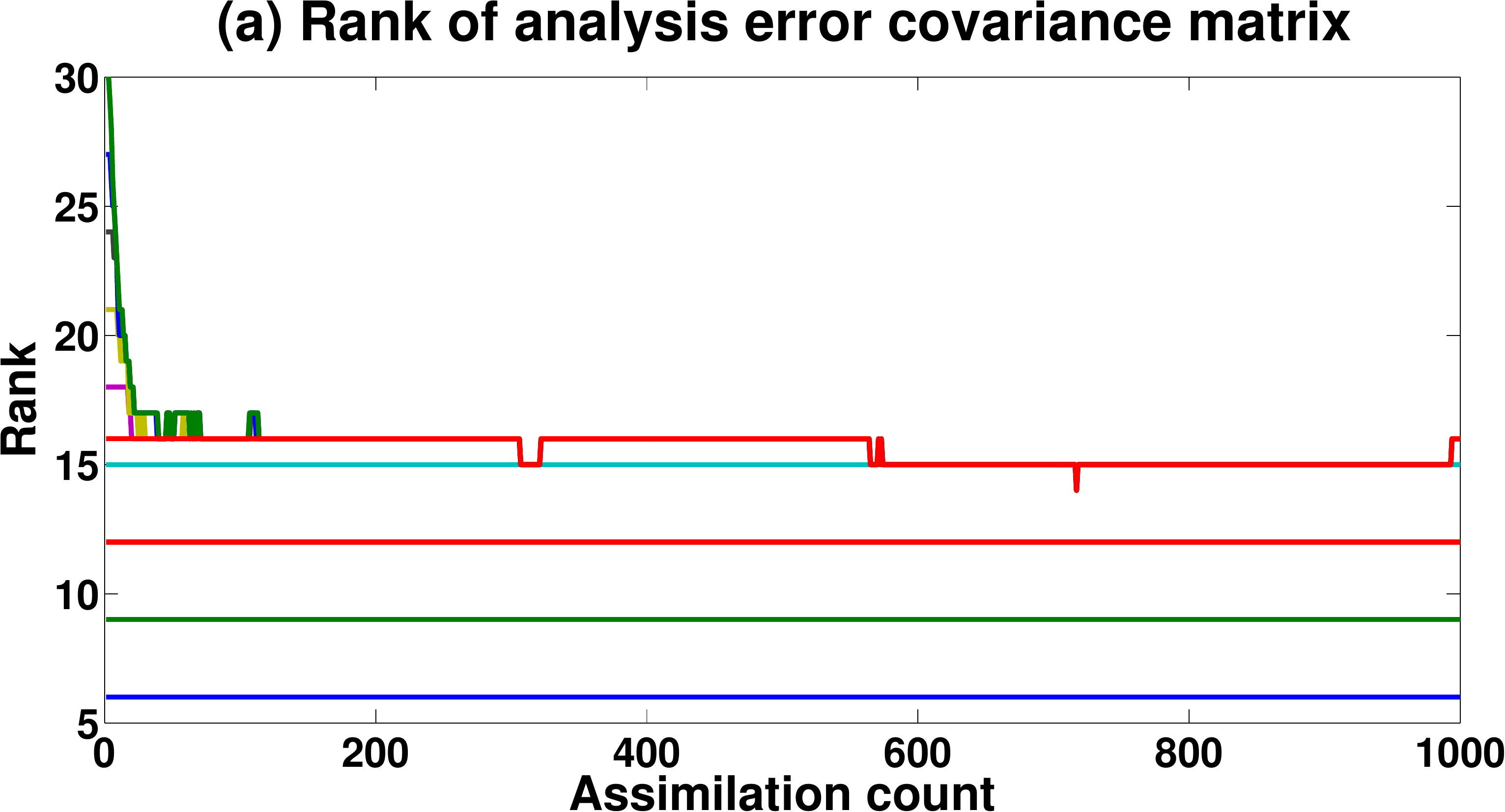}
    &
    \includegraphics[width=0.46\textwidth]{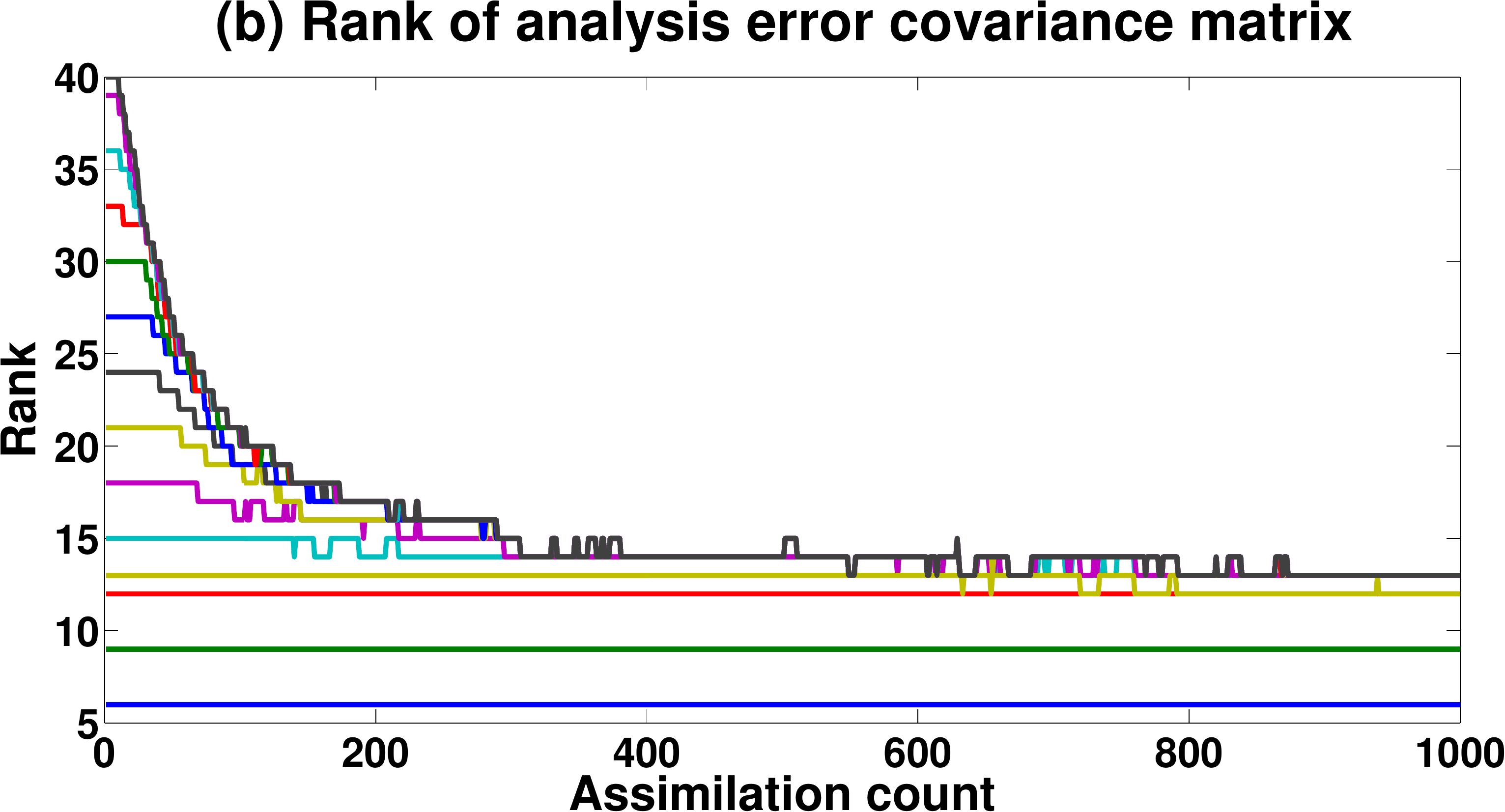}
  \end{tabular}
  \caption{\label{fig:1rankvsn} Rank of $\bP^{\rm a}_k$ as a function of $k$ for several choices of the rank $r_0$ of
    $\bP^{\rm a}_0$ (various colors) for two systems, one with random propagators $\bM_k$ (a) with $n_0 = 16$ and
    another with propagators which are a linearization of Lorenz-95 around a trajectory on the attractor (b) with $n_0 =
    14$.}
\end{figure*}

Figure~\ref{fig:1rankvsn} plots the rank of $\bP^{\rm a}_k$ as a function of $k$, where various choices of the rank
$r_0$ of $\bP^{\rm a}_0$ are shown by various colors in the figure. Note that we actually plot the number of eigenvalues
greater than a threshold of $10^{-10}$.

Panel (a) of Fig.~\ref{fig:1rankvsn} shows the case of random propagators (Exp2) $\bM_k$, which has $n=30$ and $n_0=16$,
i.e., the number of non-negative Lyapunov exponents is $16$.  Panel (b) refers to the case Exp3 of linearization of
Lorenz-95 with $F=8$ around a trajectory on its attractor with $n=40$ and $n_0 = 14$.  We see that if $r_0 < n_0$, then
the rank of $\bP^{\rm a}_k$ is constant and equal to the initial rank $r_0$. On the other hand, if $r_0 \ge n_0$, then
$r_0 - n_0$ eigenvalues values approach zero, $n_0 - 1$ eigenvalues remain non-zero, while one eigenvalue fluctuates and
it is unclear whether it will approach zero or indeed remain non-zero. It very likely corresponds to the neutral
direction along which convergence of $\bP_k$ is very slow even if well observed, as discussed in section
\ref{sec:asymptote}.

\begin{figure*}[htbp]
  \begin{tabular}{cc}
    \includegraphics[width=0.46\textwidth, clip=True]{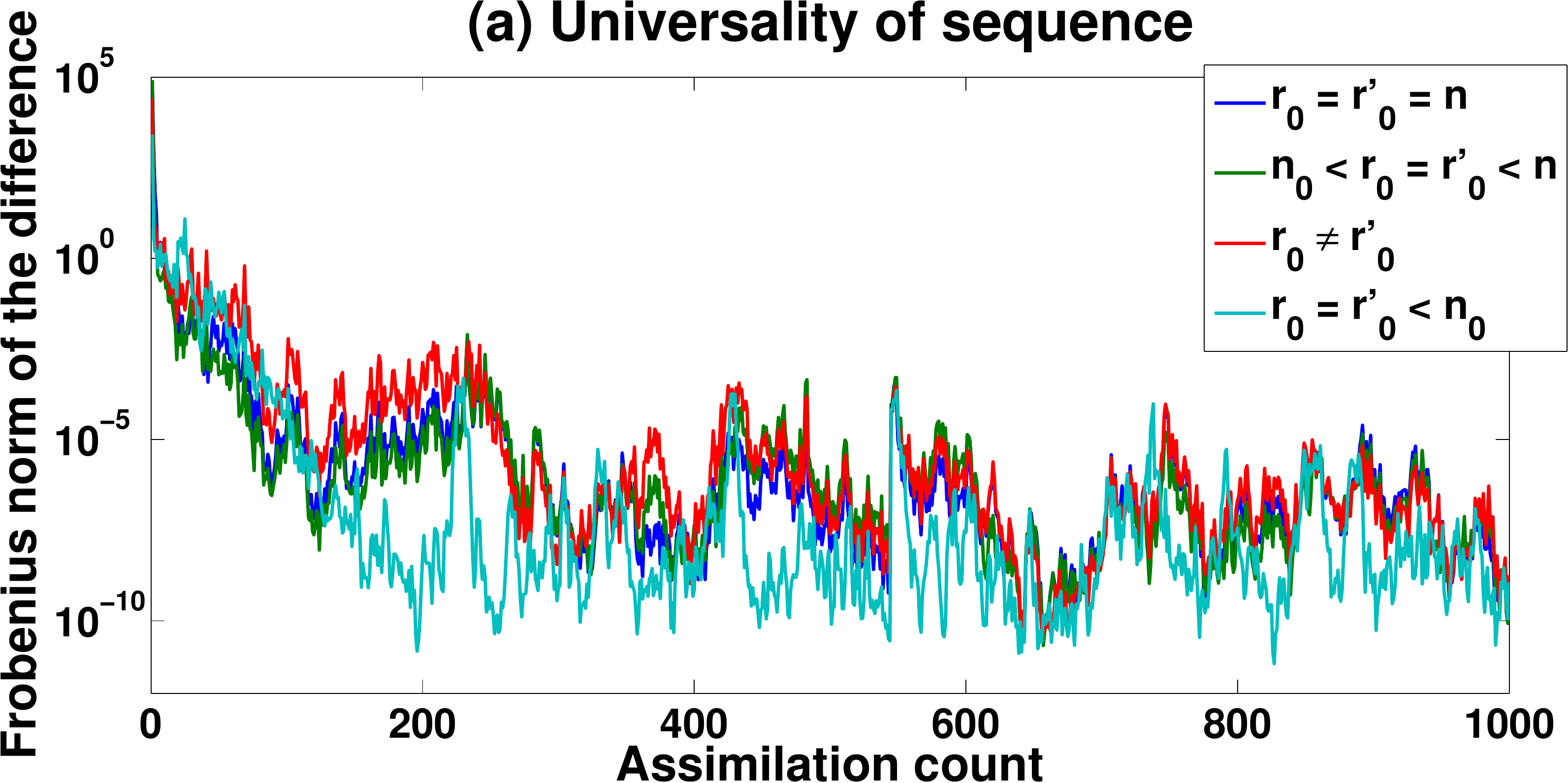}
    &
    \includegraphics[width=0.46\textwidth, clip=True]{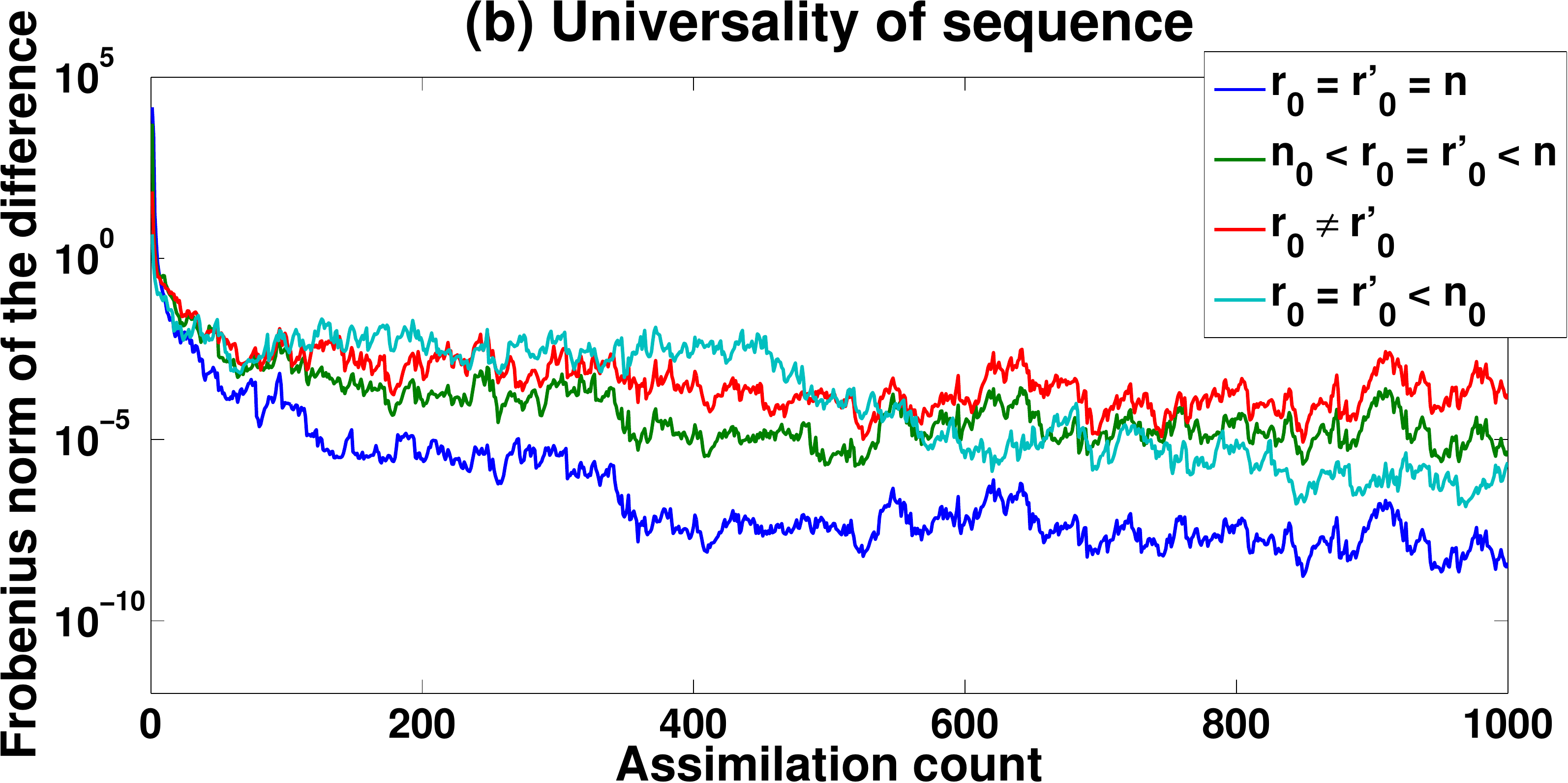}
  \end{tabular}
  \caption{\label{fig:2diff} Frobenius norm of the difference, i.e., $\|\bP^{\rm a}_k - \bP^{'\rm a}_k\|$ for two
    sequences of analysis covariances matrices starting with different initial conditions $\bP^{\rm a}_0$ and
    $\bP^{'\rm a}_0$ for the case of random propagators (a) with $n_0 = 16$ and Lorenz-95 linearization (b) with $n_0
    = 14$.}
\end{figure*}

In the next numerical experiment, we generate two sequences of analysis covariances $\bP^{\rm a}_k$ and $\bP^{'\rm a}_k$
starting from two different initial conditions $\bP^{\rm a}_0$ and $\bP^{'\rm a}_0$, respectively.
Figure~\ref{fig:2diff} shows the Frobenius norm of the difference between analysis covariances, i.e., $\|\bP^{\rm
  a}_k - \bP^{'\rm a}_k\|$, as a function of $k$.
Four cases are considered in Fig.~\ref{fig:2diff}:
\begin{enumerate}
\item $r_0 = r_0' = n$ when the initial ranks are the same and equal to the state dimension (blue line);
\item $n_0 < r_0 = r_0' < n$ when the initial covariance matrices are rank-deficient with same ranks greater than $n_0$ (green
  line);
\item $r_0 \ne r_0'$ and $n_0 < r_0 , r_0' < n$ when the initial ranks are unequal but both ranks are greater than $n_0$
  (red line);
\item $r_0 = r_0' < n_0$ when the initial ranks are the same and less than $n_0$ (teal line).
\end{enumerate}
In all these cases, we see that the norm of the difference approaches zero within the numerical accuracy, fluctuating
between $10^{-8}$ and $10^{-3}$, i.e., for large $k$, $\bP^{\rm a}_k \approx \bP^{'\rm a}_k$.  Thus the sequence
$\bP^{\rm a}_k$ is equivalent to a sequence of covariance matrices all of rank $s = \min\{r_0, n_0\}$, independent of
the initial condition $\bP^{\rm a}_0$, but of course dependent on the dynamics $\bM_{k}$, the observations $\bH_{k}$ and
their error covariances $\bR_{k}$.

The asymptotic covariance matrices are most easily represented in the basis of the BLVs. As proven mathematically in
section \ref{sec:collapse} in the case of strong collapse which occurs here because the systems are sufficiently
observed, these covariance matrices have column spaces corresponding to the span of the most unstable BLVs and their
null space subsumes the span of the stable BLVs.  This can be seen by looking at the projection of these covariance
matrices $\bP^{\rm a}_k$ onto the BLVs $\bu^k_1,\cdots,\bu^k_n$ at time $t_k$.

\begin{figure*}[htbp]
  \begin{tabular}{cc}
    \includegraphics[width=0.46\textwidth]{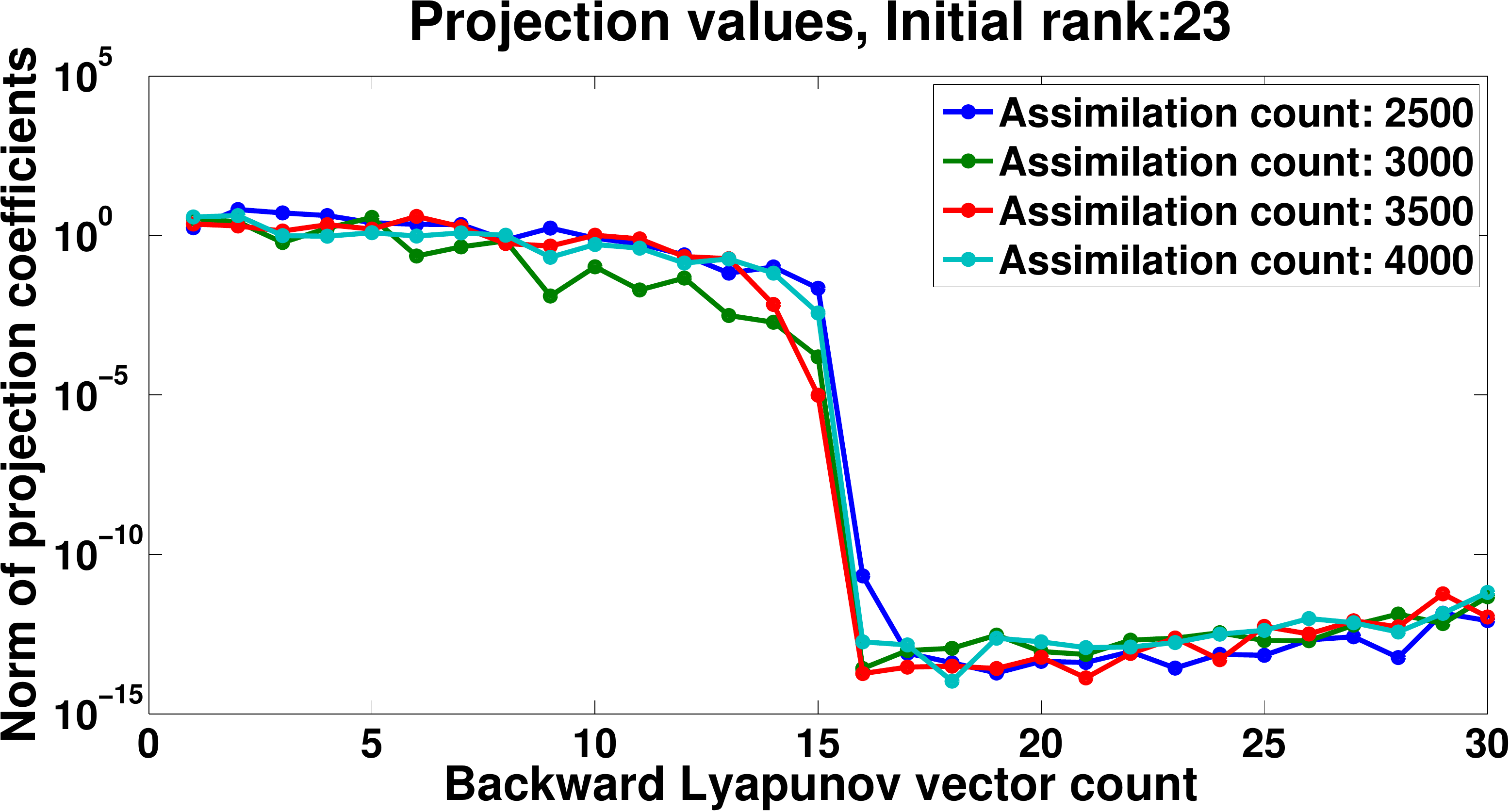}
    & \includegraphics[width=0.46\textwidth]{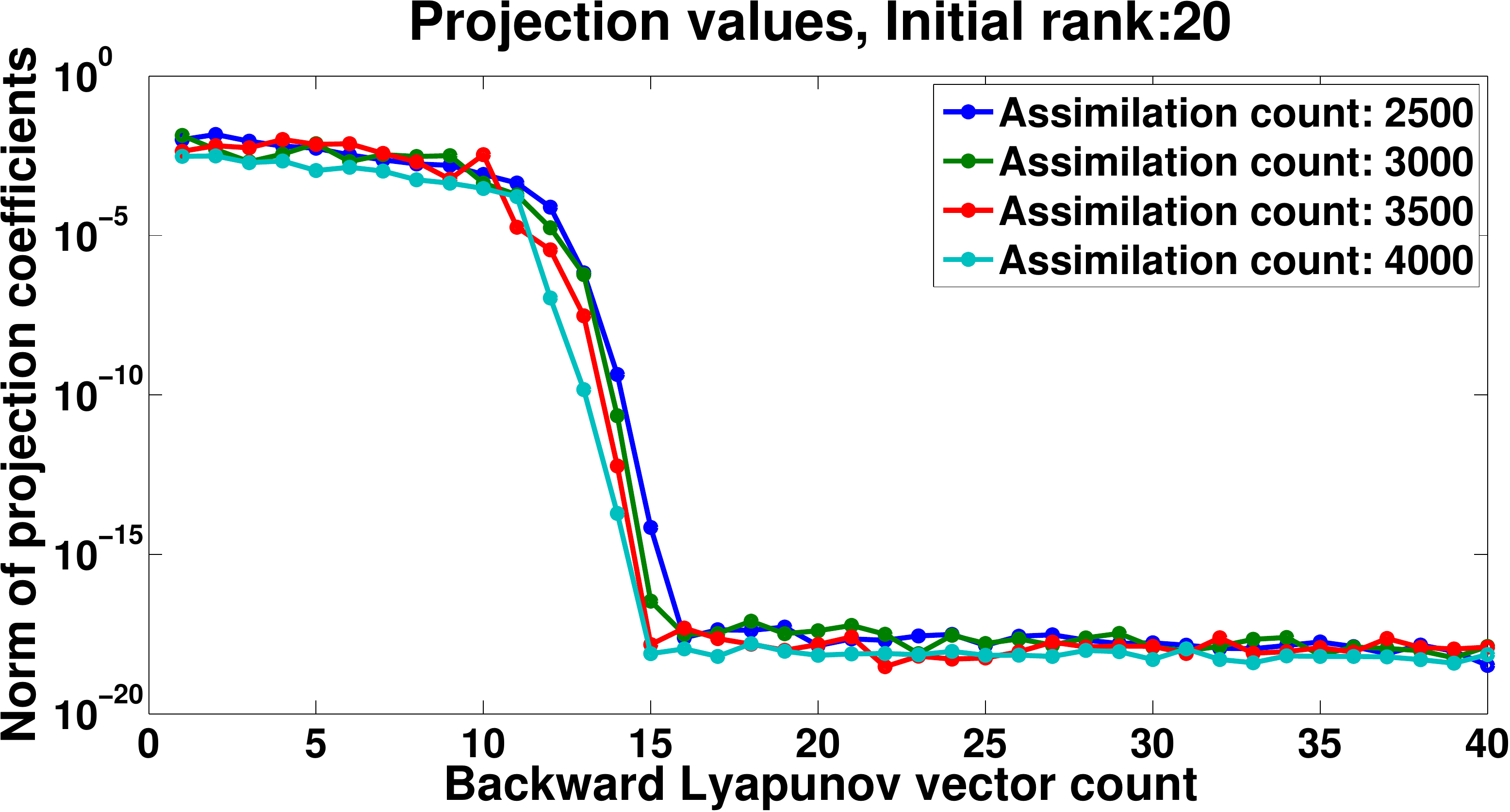} \\
    \includegraphics[width=0.46\textwidth]{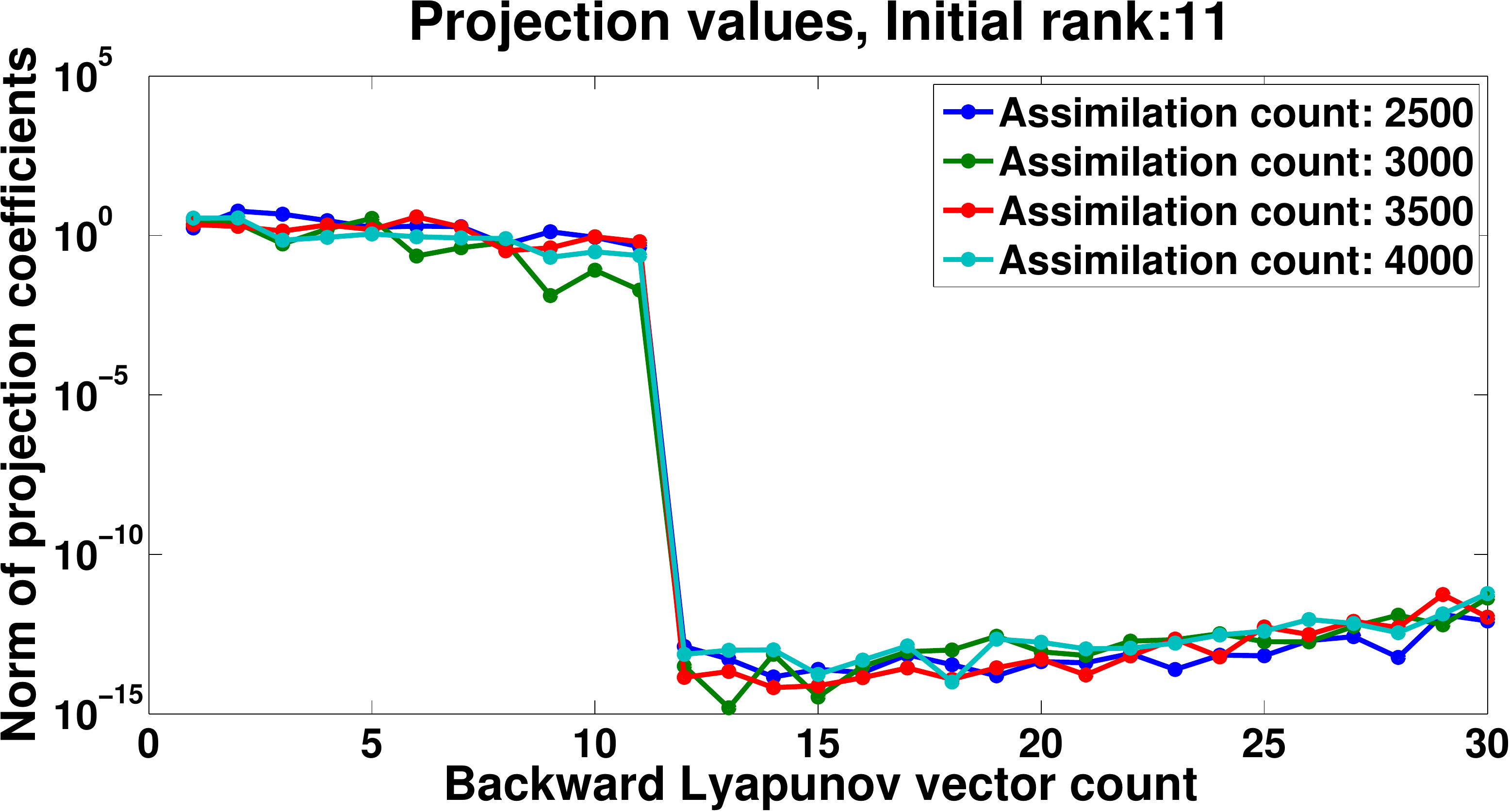}
    & \includegraphics[width=0.46\textwidth]{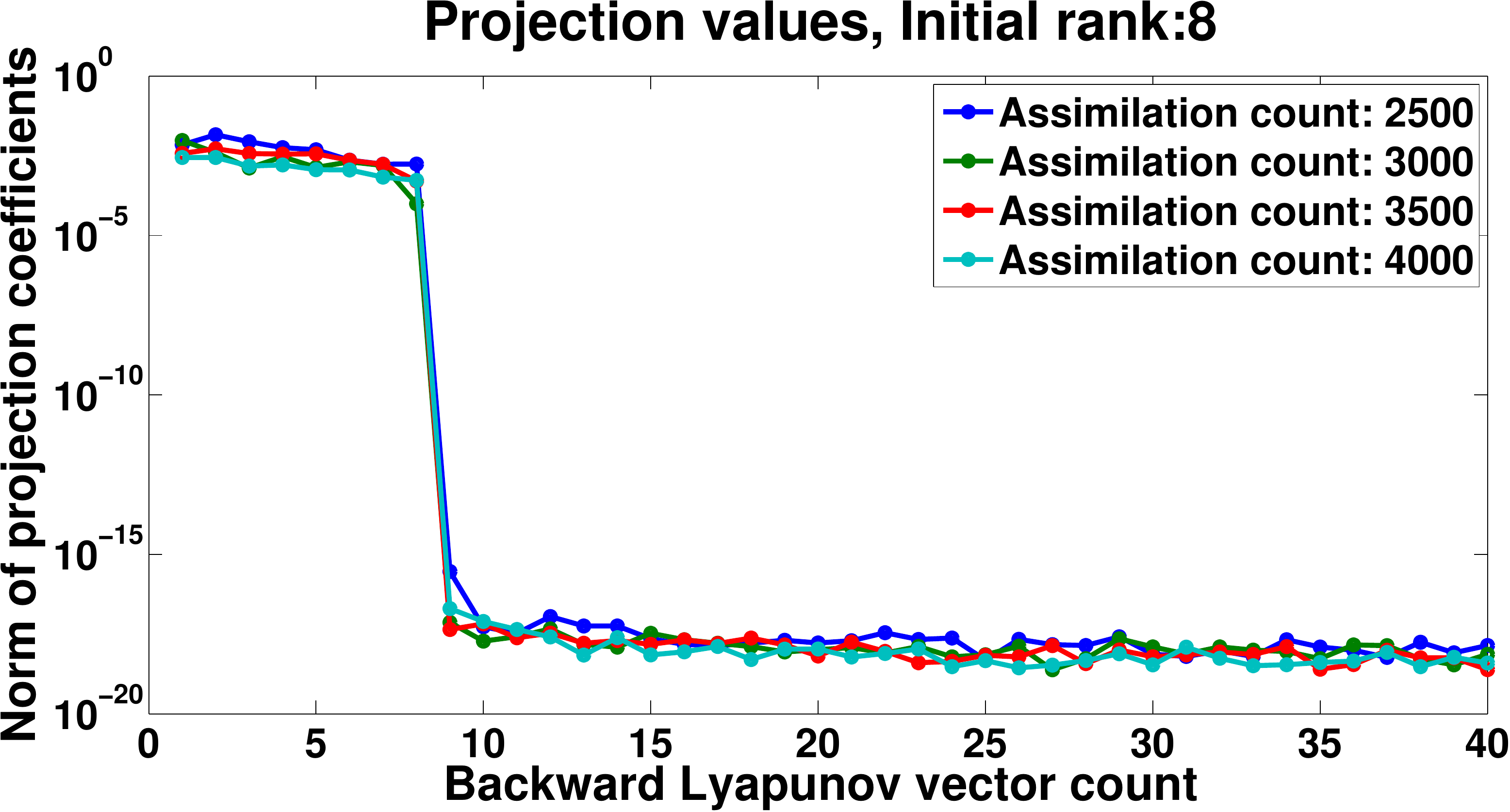} \\
    \includegraphics[width=0.46\textwidth]{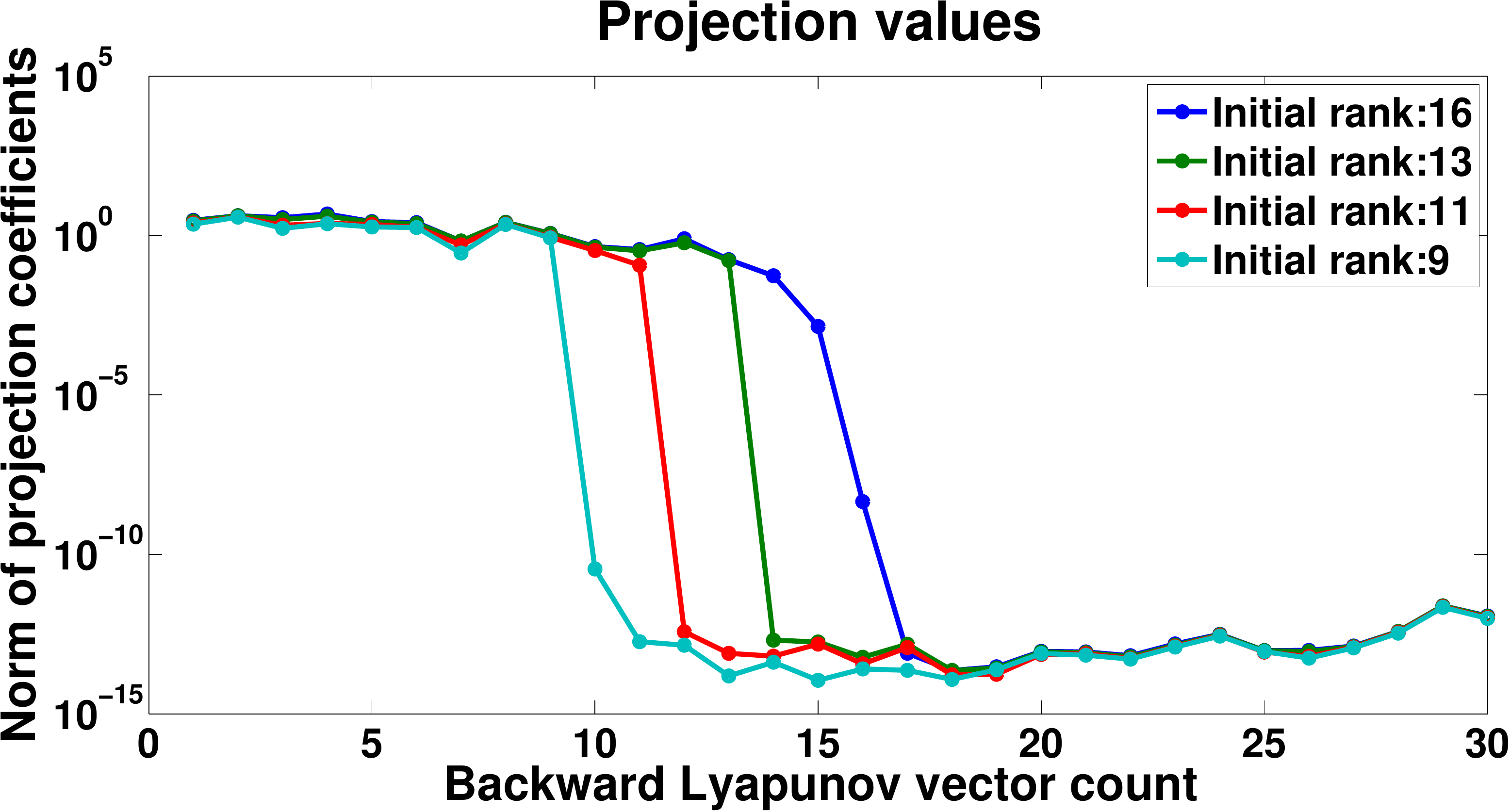} 
    & \includegraphics[width=0.46\textwidth]{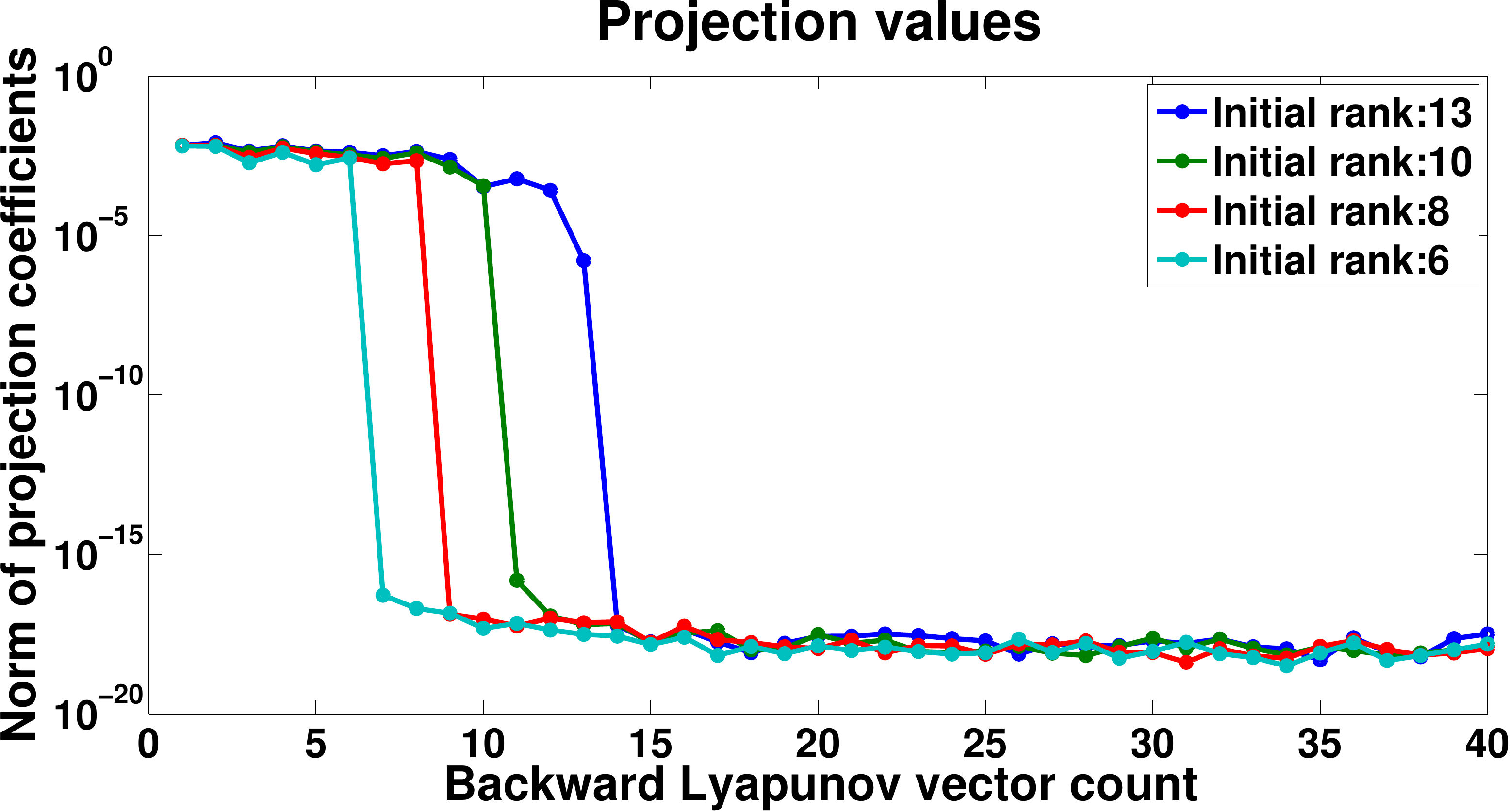} \\
  \end{tabular}
  \caption{\label{fig:3proj}
    Projections of covariance matrices $\bP^{\rm a}_k$ onto the BLVs $\bu^k_1,\cdots,\bu^k_n$ for system with random
    propagators (left column, $n = 30, n_0=16$) and linearization of Lorenz-95 (right column, $n = 40, n_0=14$).}
\end{figure*}

Figure~\ref{fig:3proj} shows these projections for four different values of $k = 2500$, $3000$, $3500$, $4000$ for the
cases $r_0 \ge n_0$ (top row) and $r_0 < n_0$ (middle row).  The Exp2 and Exp3 cases are displayed in the left and right
column panels, respectively.  Note that the Lyapunov vectors are ordered from the largest to the smallest Lyapunov
exponents.  This is also clearly seen from the bottom row of the same Fig.~\ref{fig:3proj} which shows these projections
at a fixed time $k = 5000$ for various initial ranks $r_0$ which are equal to or less than $n_0$.

\subsection{Low-rank asymptotic covariance for autonomous systems}
\label{sec:auto}

\begin{figure*}[htbp]
  \begin{tabular}{cc}
    \includegraphics[width=0.46\textwidth]{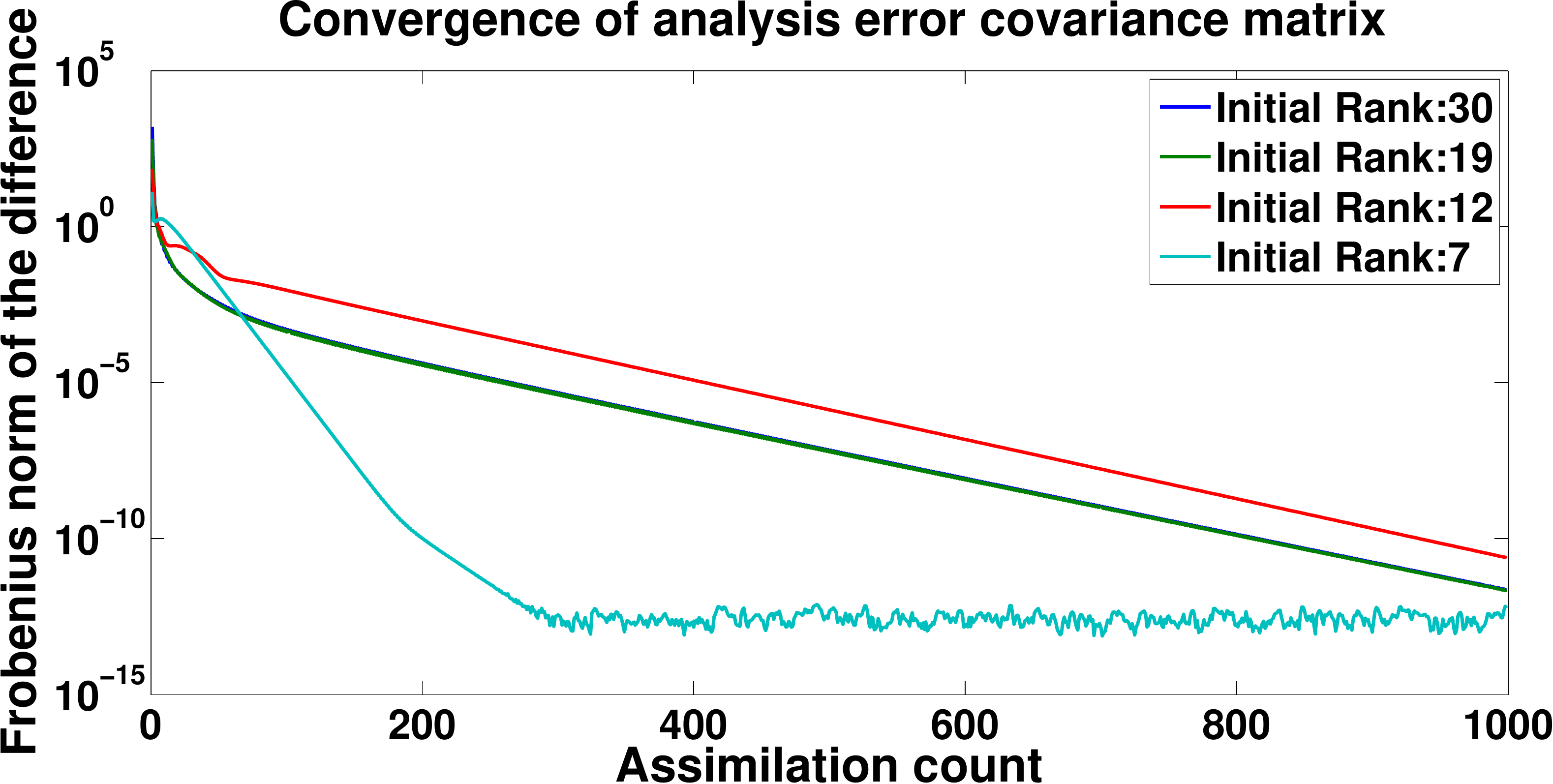}
    &
    \includegraphics[width=0.46\textwidth]{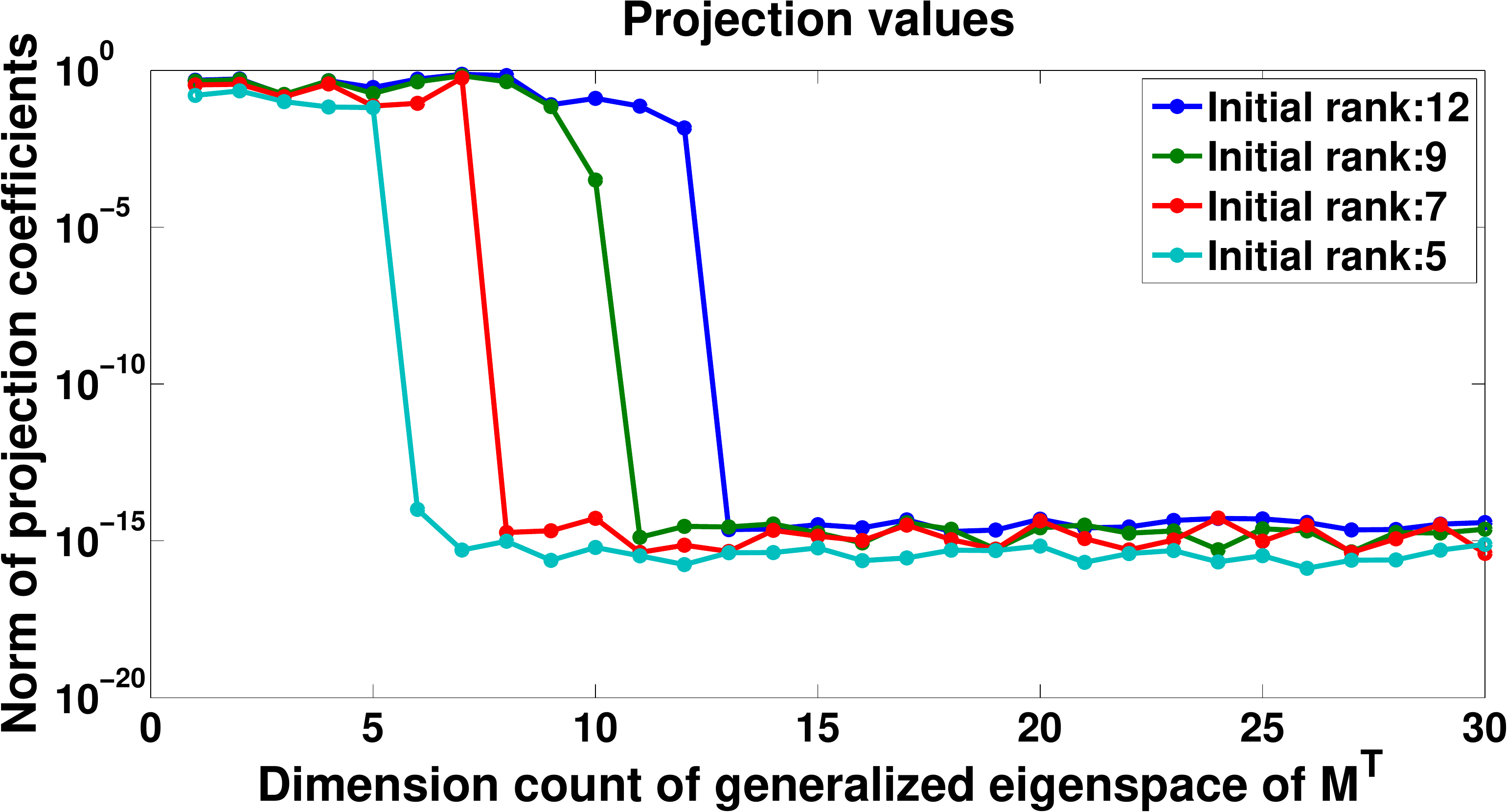}
  \end{tabular}
  \caption{\label{fig:4auto} Frobenius norm of the consecutive difference, i.e., $\|\bP^{\rm a}_k - \bP^{\rm
      a}_{k-1}\|$ for several choices of rank $r_0$ of the initial condition $\bP^{\rm a}_0$ (left panel) and
    projections onto the generalized eigenspace of $\bM^\T$ (right panel) for the autonomous system with $n = 30, n_0
    = 13$.}
\end{figure*}

The last set of numerical results illustrates the asymptotic convergence of the analysis covariances for the case of
autonomous systems.  The results are very similar to those of the non-autonomous systems and a summary is presented in
Fig.~\ref{fig:4auto}.  The left panel shows the Frobenius norm of the difference $\bP^{\rm a}_{k+1} - \bP^{\rm a}_k$ of
the analysis covariance matrices at consecutive time instances. The figure clearly shows this difference going to zero
and thus by Cauchy's convergence criterion, the sequence of the analysis covariance matrices converges.  Different lines
are meant for cases of different initial ranks. The right panel shows the projections onto the BLVs which also span the
generalized eigenspace of $\bM^\T$ \cite{gurumoorthy2017}, for four cases with different initial rank $r_0$, and these
results are very similar to those shown in the bottom row of Fig.~\ref{fig:3proj}.

\section{Conclusion}
\label{sec:concl}

We have shown that, for perfect linear dynamics and observation operator, and for any initial error covariance matrix,
the solution of the KF covariance equation converges onto the unstable-neutral subspace of the dynamics. The
rate of such convergence has also been provided. Moreover, we have shown that under reasonable assumptions there exists
a universal sequence, independent of the initial condition, toward which the Kalman filter error covariance converges if
the system is sufficiently observed and if the column space of the initial error covariance has a non-zero projection on
all the unstable and neutral FLVs.  These results were obtained after proving an analytical
expression of the covariances at any time in terms of the initial covariances.  Numerical experiments were used to
further corroborate and illustrate the mathematical statements.  These results complete and generalize those in
\cite{gurumoorthy2017} and altogether lay the mathematical foundation of the methods that rely on the assimilation in
the unstable subspace \cite{palatella2013a}.

It should be stressed in this conclusion that we have also obtained alternative square root formulas for
(\ref{eq:final2},\ref{eq:final},\ref{eq:bound1},\ref{eq:aseq}), which means that those are written in terms of
$\bX_k$ the square root factor of $\bP_k=\bX_k\bX_k^\T$ rather than $\bP_k$.  Those would have made an even stronger
connection with ensemble filters since the columns of $\bX_k$ can be seen as state perturbations associated with an
ensemble of state vectors. However, the derivations that use the square root factors turned out to be equivalent or
longer compared to working on $\bP_k$ directly. They do not bring in new insights for the purposes of this paper
compared to the derivations presented here. That is why we did not opt for the square root approach in this
paper. Square root generalizations of these formulas will be introduced elsewhere \cite{bocquet2017} in a nonlinear
context where they are more relevant.

Yet, this work leaves unresolved some key issues that are worth investigating in the perspective of the design of
reduced-order algorithms applicable to practical situations.  Specific lines of development include the treatment of
model error and the extension to nonlinear dynamics. This latter problem stimulates indeed an intriguing, albeit
necessary, direction of study whose main difficulty stands on the fact that the unstable subspace is, in this nonlinear
case, no longer globally defined but a function of the underlying trajectory. Both lines of research may lead to
interesting methodological and mathematical developments and are central in DA.

In our view the present results are also relevant to the field of ensemble-based DA algorithms for the geosciences, or,
more generally, to the uncertainty quantification and DA methods in complex high-dimensional and big data problems, of
which DA for the geosciences is a prototypical example. We believe so for two distinctive reasons. First the present
findings on the error covariance projection onto the unstable-neutral subspace provide a natural rationale to interpret
a stream of numerical evidence that relates the minimum ensemble size to achieve a satisfactorily estimate of the
system's state, with the number of unstable directions of the underlying dynamics
\cite{carrassi2009,ng2011,bocquet2014}.  Second, this study encourages a research effort toward EnKF formulations that
incorporate the information on the unstable subspace explicitly in the design and choice of the ensemble, possibly in
combination with localization techniques widely used to artificially increase the rank of the ensemble-based error
covariance matrices.

While a specific recipe for a formulation of the EnKF under this framework is still part of the authors’ ongoing
research, some preliminary considerations can nevertheless be put forward. First, the convergence onto the unstable
subspace for EnKF covariance can only be obtained for the class of the deterministic EnKFs, as confirmed by the
numerical results performed by the authors \cite{bocquet2017}. Second, for the sake of the feasibility, the explicit
use of the unstable subspace in the filter design for high-dimensional applications, must necessarily rely on an
efficient computation of such a subspace. Recent developments along this line \cite{wolfe2007,ginelli2013} thus appear
favorable and further support our current research.

It is finally worth mentioning another appealing research direction: the extension of the present framework to fully
Bayesian DA methods typically preferable in the presence of strong nonlinearities and/or non-Gaussian error
\cite{bocquet2010}. Besides the aforesaid difficulty inherent to nonlinear dynamics, the additional problem here is on
how to link the geometrical (in the phase-space) features of the unstable subspace to the conditioning of a pdf, that is
the generalization to the fully Bayesian framework, of projecting the error covariances onto the unstable subspace.

\appendix

\section{Deriving $\bP_k$ using the symplectic symmetry}
\label{app:symplectic}

This appendix gives an account of the linear representation of the recurrence \eqref{eq:recurrence}, which had
initially been developed as a way to solve the Riccati equation in the autonomous case \cite[and references
  therein]{anderson1979}.  We use it to give an alternative derivation of \eqref{eq:final2} and to discuss in more
detail the analytic expression in the autonomous case.  The underlying symplectic structure of the KF has
been, for instance, explored in \cite{bougerol1993, wojtowski2007}.

\subsection{General linear representation using symplectic matrices}

Let us rewrite the recurrence \eqref{eq:recurrence}:
\begin{align}
  \bP_{k} =& \bM_{k}\(\I{n}+\bP_{k-1} \bOmega_{k-1}\)^{-1}\bP_{k-1}\bM_{k}^{\T}+\bQ_k \nn
  =& \bM_{k}\bP_{k-1}\(\I{n}+\bOmega_{k-1}\bP_{k-1}\)^{-1}\bM_{k}^{\T}+\bQ_k \nn
  =& \bM_{k}\bP_{k-1}\(\bM_{k}^{-{\T}}+\bM_{k}^{-{\T}} \bOmega_{k-1}\bP_{k-1}\)^{-1} +\bQ_k \nn
  =& \(\bM_{k}\bP_{k-1} + \bQ_k\left\{\bM_{k}^{-{\T}}+\bM_{k}^{-{\T}} \bOmega_{k-1}\bP_{k-1}\right\}\)
  \(\bM_{k}^{-{\T}}+\bM_{k}^{-{\T}} \bOmega_{k-1}\bP_{k-1}\)^{-1} \nn
  =& \(\left\{\bM_{k}+\bQ_k\bM_{k}^{-{\T}} \bOmega_{k-1}\right\}\bP_{k-1} + \bQ_k\bM_{k}^{-\T}\)
  \(\bM_{k}^{-{\T}} \bOmega_{k-1}\bP_{k-1}+ \bM_{k}^{-{\T}}\)^{-1} \nn
   \triangleq & \(\bA_k \bP_{k-1} + \bB_k\)\(\bC_k\bP_{k-1}+\bD_k\)^{-1} ,
\end{align}
where we used the matrix shift lemma from the first to the second line,
and we defined block matrices $\bA_k$, $\bB_k$, $\bC_k$, $\bD_k$ in the fourth line.
Let us define \cite{anderson1979}
\be
\bZ_k \triangleq
\begin{pmatrix}
\bA_k & \bB_k \\  
\bC_k & \bD_k
\end{pmatrix}
=
\begin{pmatrix}
\bM_{k}+ \bQ_k\bM_k^{-{\T}}\bOmega_{k-1} &  \bQ_k\bM_k^{-{\T}}\\  
\bM_k^{-{\T}}\bOmega_{k-1} & \bM_k^{-{\T}},
\end{pmatrix}.
\ee
which is valid in the presence of model noise.  This matrix belongs to the symplectic group ${\rm Sp}(2n,{\mathbb R})$
since $\bZ^{-1}_k=-\bJ \bZ_k^{\T}\bJ$, where $\bJ = \begin{pmatrix} \bzero & \I{n} \\ -\I{n} & \bzero \end{pmatrix}$.  It
has a simple expression in the perfect model case:
\be
\label{eq:defZk}
\bZ_k \triangleq
\begin{pmatrix}
\bA_k & \bB_k \\  
\bC_k & \bD_k
\end{pmatrix}
=
\begin{pmatrix}
\bM_{k} & \bzero \\  
\bM_k^{-{\T}}\bOmega_{k-1} & \bM_k^{-{\T}}
\end{pmatrix}.
\ee

Furthermore, let us introduce the following matrix in $\R{2n}{n}$:
\be
\bW_k = \begin{pmatrix} \bX_k \\ \bY_k \end{pmatrix} ,
\ee
where $\bY_k$ is assumed to be invertible, which can and will be checked a posteriori, and we define the ratio $\bomega_k
= \bX_k \bY_k^{-1}$ in $\R{n}{n}$.  The $\bW_k$ are related by the defining recurrence
\be
\label{eq:symrec}
\bW_{k+1} \triangleq \bZ_k\bW_k .
\ee
We explicitly have
\be
\begin{pmatrix} \bX_{k+1} \\ \bY_{k+1} \end{pmatrix} \triangleq
\bZ_k\bW_k = 
\begin{pmatrix}
\bA_k & \bB_k \\  
\bC_k & \bD_k
\end{pmatrix}
\begin{pmatrix} \bX_k \\ \bY_k \end{pmatrix}
=
\begin{pmatrix} \bA_k\bX_k +\bB_k\bY_k \\ \bC_k\bX_k + \bD_k\bY_k \end{pmatrix},
\ee
from which it is possible to infer the following recurrence on $\bomega_k$:
\begin{align}
  \bomega_{k+1} &= \bX_{k+1} \bY_{k+1}^{-1} = (\bA_k\bX_k +\bB_k\bY_k)(\bC_k\bX_k + \bD_k\bY_k)^{-1} \nn
&= (\bA_k\bX_k\bY_k^{-1} +\bB_k)(\bC_k\bX_k\bY_k^{-1} + \bD_k)^{-1} \nn
& = (\bA_k\bomega_k +\bB_k)(\bC_k\bomega_k + \bD_k)^{-1}.
\end{align}
Hence, we can represent the nonlinear update of $\bomega_k$ by the linear recurrence
\eqref{eq:symrec}.

Now, we choose
\be
\bX_0 = \bP_0 \quad {\rm and} \quad \bY_0 = \I{n}
\ee
in order to have $\bomega_k = \bP_k$ for all $k \ge 0$, which implies that the nonlinear recurrence on $\bP_k$ can be
represented by the linear recurrence \eqref{eq:symrec}.

Insofar, no assumption on the rank of $\bP_k$ was required and, even in the presence of model noise, the linear
representation implies that $\bP_k$ has the following dependence on $\bP_0$:
\be
\bP_k = (\bA^{(k)}\bP_0 +\bB^{(k)})(\bC^{(k)}\bP_0 + \bD^{(k)})^{-1}
\ee
where the $\bA^{(k)}$, $\bB^{(k)}$, $\bC^{(k)}$, $\bD^{(k)}$ only depend on $\bOmega_l$, $\bQ_l$, and $\bM_l$, $1\le l \le k$.  Our
purpose now is to compute $\bP_k$ for any $t_k$ {\em in the perfect model} case using the linear representation
\eqref{eq:defZk}. This is the focus of the rest of this appendix.

\subsection{Solution in the autonomous case}

We consider first the autonomous case, where $\bM_k$, $\bOmega_k$ and $\bZ_k$ are all independent of time, and we can
suppress the time index from the notation.  Hence, we would compute the power iterates $\bZ^k$ of $\bZ$ (not to be
confused with the $\bZ_k$ defined in \eqref{eq:defZk}. Let us assume that $\bZ^k$ has the form
\be
\bZ^k \triangleq 
\begin{pmatrix}
\bM^k & \bzero \\  
\(\bM^{k}\)^{-{\T}}\bTheta'_k & \(\bM^{k}\)^{-{\T}}
\end{pmatrix}
.
\ee
Note that we want $\bZ^0 = \begin{pmatrix}
\I{n} & \bzero \\  
\bzero & \I{n} \end{pmatrix}$, so that $\bTheta'_0 = \bzero$. 
Then the recurrence on $\bZ^k$ imposes the recurrence on the $\bTheta'_k$
\be
\label{eq:rectheta}
\bTheta'_{k+1} = \bM^{\T}\bTheta'_k\bM + \bOmega ,
\ee
which identifies $\bTheta'_k$ with $\bTheta_k$ as defined by \eqref{eq:thetadef}.  Because $\bTheta_0$ and $\bOmega$
are symmetric, all $\bTheta_k$ for $k\ge 1$ are also symmetric.  We can see it as an arithmetico-geometric recurrence
and to solve it define by $\bPsi$ the solution of
\be
\label{eq:lyapunov}
  \bPsi = \bM^{\T} \bPsi \bM + \bOmega .
\ee 
This is the so-called discrete algebraic Lyapunov equation.
Because a solution of this equation does not always exist \cite{gajic1995}, we consider instead the recurrence
  \be
\label{eq:rectheta2}
  e^{i\varepsilon}\bTheta^\varepsilon_{k+1} = \bM^{\T}\bTheta^\varepsilon_k\bM + \bOmega,
\ee
where $0 < \varepsilon < 2\pi$. By continuity, $\bTheta_k = \lim_{\varepsilon \rightarrow 0^+} \bTheta^\varepsilon_k$.
The corresponding Lyapunov equation is
\be
\label{eq:lyapunov2}
e^{i\varepsilon} \bPsi_\varepsilon = \bM^{\T}\bPsi_\varepsilon\bM + \bOmega .
\ee
Because it is formally equivalent to
\be
\label{eq:lyapunov3}
\(e^{i\varepsilon}\I{n} - \bM^\T\otimes\bM^\T \) \mathrm{vec}\(\bPsi\) = \mathrm{vec}\(\bOmega\),
\ee
where $\mathrm{vec}\(\bPsi\)$ is the vector made from the stacked columns of $\bPsi$. Since \be
\mathrm{det} \(e^{i\varepsilon}\I{n} - \bM^\T\otimes\bM^\T\) \neq 0 ,
\ee
for any real $\bM$ and $0 < \varepsilon < 2\pi$, there exists a unique solution $\bPsi_\varepsilon$ of \eqref{eq:lyapunov3} in
$\mathbb C$.  Then, we obtain 
\be
\bTheta^\varepsilon_k = \bPsi_\varepsilon - (\bM^k)^{\T}\bPsi_\varepsilon\bM^k ,
\ee
by subtracting \eqref{eq:lyapunov2} from \eqref{eq:rectheta2} and then iterating.
As a consequence, the following construction of a solution is always valid:
\be
\bTheta_k = \lim_{\varepsilon \rightarrow 0^+} \left\{ \bPsi_\varepsilon - e^{-ik\varepsilon}\(\bM^k\)^\T\bPsi_\varepsilon\bM^k \right\} \, .
\ee
Hence,
\be
\bZ^k =
\begin{pmatrix}
\bM^k & \bzero \\  
\(\bM^{k}\)^{-{\T}}\bTheta_k & \(\bM^{k}\)^{-{\T}}
\end{pmatrix}
.
\ee
Using the linear representation leads to
\begin{align}
\begin{pmatrix} \bX_k \\ \bY_k \end{pmatrix}
&= \begin{pmatrix}
\bM^k & \bzero \\  
\(\bM^{k}\)^{-{\T}}\bTheta_k & \(\bM^{k}\)^{-{\T}}
\end{pmatrix}
\begin{pmatrix} \bP_0 \\ \I{n} \end{pmatrix} \nn
&=
\begin{pmatrix} \bM^k\bP_0 \\ \( \bM^{k}\)^{-{\T}}\bTheta_k \bP_0 + \(\bM^{k}\)^{-{\T}} \end{pmatrix} .
\end{align}
Using $\bP_k = \bX_k\bY_k^{-1}$, we conclude
\be
  \bP_k = \bM^k\bP_0 \left[ \bTheta_k \bP_0 + \I{n} \right]^{-1}\(\bM^{k}\)^{\T} .
\ee

\subsection{Solution in the non-autonomous case}

In the non-autonomous case, we need to define
\be
\bZ^{(k)} \triangleq \bZ_k \bZ_{k-1}\cdots \bZ_0 .
\ee
The product is still in the symplectic group and of the form
\be
\bZ^{(k)} \triangleq 
\begin{pmatrix}
\bM_{k:0} & \bzero \\  
\bGamma'_k \bM_{k:0} & \bM_{k:0}^{-\T}
\end{pmatrix},
\ee
which leads to the following recurrence on $\bGamma'_k$:
\be
\bGamma'_{k+1} = \bM^{-\T}_{k+1}\(\bGamma'_k + \bOmega_k \) \bM^{-1}_{k+1}.
\ee
The finite-time solution to this recurrence is
\be
\bGamma'_k = \sum_{l=0}^{k-1} \bM_{k:l}^{-\T} \bOmega_l \bM_{k:l}^{-1}
\ee
which coincides with the definition of $\bGamma_k$ in \eqref{eq:gammadef}.  Hence, we have an expression for
$\bZ^{(k)}$. We can use it to obtain a solution for the recurrence on $\bP_k$ using the linear representation
\be
\bZ^{(k)}\begin{pmatrix} \bP_0 \\ \I{n} \end{pmatrix} =
\begin{pmatrix} \bM_{k:0}\bP_0 \\ \bGamma_k \bM_{k:0}\bP_0 + \bM_{k:0}^{-\T} \end{pmatrix}
\ee
from which we obtain
\begin{align}
\label{eq:na}
\bP_k &= \bM_{k:0}\bP_0 \left[ \bGamma_k \bM_{k:0}\bP_0 + \bM_{k:0}^{-\T}\right]^{-1} \nn
&= \bM_{k:0}\bP_0  \bM_{k:0}^{\T} \left[ \I{n} + \bGamma_k \bM_{k:0}\bP_0 \bM_{k:0}^{\T}\right]^{-1} \nn
&= \bM_{k:0}\bP_0  \left[ \I{n} + \bM_{k:0}^{\T} \bGamma_k \bM_{k:0}\bP_0 \right]^{-1}\bM_{k:0}^{\T}
\end{align}
which coincides with \eqref{eq:final2} and \eqref{eq:final3}.

\section{A few useful properties of the symmetric positive (semi-) definite matrices}
\label{app:cone}

Here we provide a selection of definitions and results about the symmetric positive (semi-)definite matrices that we use
in this paper. An introduction and detailed proofs of several of these results can be found in \cite[chapter 6]{zhang1999}.

\begin{enumerate}
\item
  
The partial ordering on $\cone{n}$ is defined by, for any $\bA$ and $\bB$ in $\cone{n}$, $\bA \le \bB$ if and only if for
all $\bx \in \Rn$, $\bx^\T\bA\bx \le \bx^\T\bB\bx$.

\item
If $\bA$ and $\bB$ are in $\cone{n}$ and $\bG$ is in $\R{q}{n}$, $q \in {\mathbb N}$, we have that $\bA \le \bB$ implies
$\bG \bA \bG^\T \le \bG \bB \bG^\T$ which is immediate from the previous definition of the partial ordering.

\item
If $\bA$ and $\bB$ are in $\conep{n}$, $\bA \le \bB$ is equivalent to $\bA^{-1} \ge \bB^{-1}$.  This can be shown using
the double diagonalization theorem which states that there exists an invertible matrix $\bG$ such that $\bG\bA\bG^\T$ and
$\bG\bB\bG^\T$ are both diagonal.

\item
If $\bA$, $\bB$, and $\bC$ are in $\cone{n}$, by $\bA \le \min \{\bB,\bC\}$ we mean that for all $\bx \in \Rn$,
$\bx^\T\bA\bx \le \min \{\bx^\T\bB\bx,\bx^\T\bC\bx\}$.

\item
  If $\bA$ is in $\cone{n}$, it has the eigendecomposition $\bA = \sum_{i=1}^n \sigma_i \bv_i\bv_i^\T$, with $\sigma_i
  \ge 0$ and $\left\{\bv_i\right\}_{1\le i \le n}$ an orthonormal basis. Let $\sigma_\mathrm{max} = \max_{1\le i \le n}
  \sigma_i $ and $\sigma_\mathrm{min} = \min_{1\le i \le n} \sigma_i $.  Any $\bx$ in $\Rn$ can decompose on the
  eigenvectors of $\bA$: $\bx = \sum_{i=1}^n \(\bv_i^\T\bx\) \bv_i$. As a consequence, one has
  \begin{align}
    \bx^\T\bA\bx = \sum_{i=1}^n \sigma_i  \(\bv_i^\T\bx\)^2
    & \le \sigma_\mathrm{max} \sum_{i=1}^n \(\bv_i^\T\bx\)^2
    =  \sigma_\mathrm{max} \bx^\T\bx \nn
    & \ge \sigma_\mathrm{min} \sum_{i=1}^n \(\bv_i^\T\bx\)^2 = \sigma_\mathrm{min} \bx^\T\bx 
  \end{align}
which leads to $\sigma_\mathrm{min} \I{n} \le \bA \le \sigma_\mathrm{max} \I{n}$. Further, as
\be
\left\| \bA \bx \right\|^2 = \sum_{i=1}^n \sigma_i ^2 \(\bv_i^\T\bx\)^2 ,
\ee
it follows that $\bx^\T\bA\bx = 0 \iff \bA \bx = \bzero$.

Now, assume $\left\{\bA_k\right\}_{k \in {\mathbb N}}$ is a uniformly bounded sequence in $\cone{n}$ and
$\left\{\bx_k\right\}_{k \in {\mathbb N}}$ is a uniformly bounded sequence in ${\mathbb R}^n$.  Then $\lim_{k
  \rightarrow \infty} \bA_k \bx_k = \bzero$ implies that $\lim_{k \rightarrow \infty} \bx_k^\T \bA_k \bx_k = 0$ by
virtue of the boundedness of $\bx_k$.  Owing to the uniform boundedness of $\bA_k$, we introduce $\sigma = \sup_{k \in
  {\mathbb N}, 1\le i\le n}\sigma_{k,i} < \infty$ and obtain
\begin{align}
  \left\| \bA_k \bx_k \right\|^2 =& \sum_{i=1}^n \sigma_{k,i} ^2 \(\bv_{k,i}^\T\bx_k\)^2  \nn
  \le&  \sigma \sum_{i=1}^n \sigma_{k,i} \(\bv_{k,i}^\T\bx_k\)^2 \le \sigma \bx_k^\T\bA_k\bx_k .
\end{align}
Hence, $\lim_{k \rightarrow \infty} \bx_k^\T \bA_k \bx_k = 0$ implies that $\lim_{k \rightarrow \infty} \bA_k \bx_k =
\bzero$.  It follows that $\lim_{k \rightarrow \infty} \bx_k^\T\bA_k\bx_k = 0 \iff \lim_{k \rightarrow \infty} \bA_k
\bx_k = \bzero$.  In particular, if the diagonal entry $[\bA_k]_{ii}$ asymptotically vanishes, the associated row
$[\bA_k]_{i\cdot}$ and column $[\bA_k]_{\cdot i}$ asymptotically vanish. Accordingly, if a given diagonal block of the
$\bA_k$ asymptotically vanishes, the off-diagonal blocks with the same row and column indices as the diagonal block
asymptotically vanish.

\item
Let $\bA \in \cone{n}$ and $\alpha \geq 0$ be a constant. If there is a subspace $\W \subseteq \Rn$ of dimension $\s \ge
1$ such that for all unit vectors $\bh \in \W$, $\bh^{\T} \bA \bh \leq \alpha$, then $\bA$ has at least $\s$ of its
eigenvalues less than or equal to $\alpha$.

To see this, decompose $\bA = \sum_{i=1}^n \sigma_i \bv_i\bv_i^\T$ in its orthonormal eigenbasis where $\sigma_i \ge 0$
and ordered as $\sigma_1 \ge \sigma_2 \ge \cdots \ge \sigma_n$. Consider $\V$ the $(s-1)$-dimensional subspace span of
$\left\{ \bv_{n-s+2}, \ldots, \bv_{n} \right\}$, which we take to be the null space if $s=1$. The orthogonal subspace
$\V^\perp$ of $\V$ is of dimension $n-s+1$. The intersection $\W \cap \V^\perp$ is of dimension at least $1$. Let us
pick $\bh$ of Euclidean norm $1$ in this intersection. We have
\begin{align}
  \alpha \ge \bh^\T\bA\bh &= \sum_{i=1}^n \sigma_i  (\bh^\T\bv_i)^2 = \sum_{i=1}^{n-s+1} \sigma_i  (\bh^\T\bv_i)^2 \nn
  &\ge \sigma_{n-s+1} \sum_{i=1}^{n-s+1} (\bh^\T\bv_i)^2 = \sigma_{n-s+1} .
\end{align}
Hence $\alpha \ge \sigma_{n-s+1} \ge \cdots \ge \sigma_n$.
 
\end{enumerate}

\section{Matrix shift lemma}
\label{app:sml}
Let $\bA \in \R{n}{m}$ and $\bB \in \R{m}{n}$.  Assuming $x \mapsto f(x)$ can be written as a formal power series,
i.e., $f(x)=\sum_{i=0}^{\infty}a_i x^i$, one has $\bA f(\bB\bA)= \sum_{i=0}^{\infty}a_i \bA(\bB\bA)^i =
\sum_{i=0}^{\infty}a_i (\bA\bB)^i\bA = f(\bA\bB)\bA$. This proves the matrix shift lemma, i.e. $\bA
f(\bB\bA)=f(\bA\bB)\bA$.  In the special case that $f(x) = (1+x)^{-1}$, this property in fact holds for any matrix
without consideration of the radius of convergence of the power series.  Assuming $(\I{m} +\bB\bA)^{-1}$ and $(\I{n} +
\bA\bB)^{-1}$ exist then
\be
\bA(\I{m} +\bB\bA)^{-1} = (\I{n} + \bA\bB)^{-1}(\bA + \bA\bB\bA)(\I{m} +\bB\bA)^{-1} = (\I{n} + \bA\bB)^{-1}\bA .
\ee

\section*{Acknowledgments}
The authors thank Sebastian Reich acting as reviewer and an anonymous reviewer for their comments and suggestions.
A. Carrassi has been partly funded by the Nordic Centre of Excellence {\em EmblA} of the Nordic Countries Research
Council, NordForsk and, along with C. Grudzien, by the project REDDA of the Norwegian Research Council (contract:
250711).  A. Apte and K. Gurumoorthy benefited from the support of the Airbus Group Corporate Foundation Chair in
Mathematics of Complex Systems established in ICTS-TIFR.

\bibliographystyle{siamplain}
\bibliography{references}

\begin{thebibliography}{10}

\bibitem{anderson1979}
{\sc B.~D.~O. Anderson and J.~B. Moore}, {\em Optimal Filtering},
  Prentice-Hall, Inc, Englewood Cliffs, New Jersey, 1979.

\bibitem{bocquet2017}
{\sc M.~Bocquet and A.~Carrassi}, {\em Four-dimensional ensemble variational
  data assimilation and the unstable subspace}, Tellus A, 69 (2017),
  p.~1304504.

\bibitem{bocquet2010}
{\sc M.~Bocquet, C.~A. Pires, and L.~Wu}, {\em Beyond {G}aussian statistical
  modeling in geophysical data assimilation}, Mon. Wea. Rev., 138 (2010),
  pp.~2997--3023.

\bibitem{bocquet2015}
{\sc M.~Bocquet, P.~N. Raanes, and A.~Hannart}, {\em Expanding the validity of
  the ensemble {K}alman filter without the intrinsic need for inflation},
  Nonlin. Processes Geophys., 22 (2015), pp.~645--662.

\bibitem{bocquet2014}
{\sc M.~Bocquet and P.~Sakov}, {\em An iterative ensemble {K}alman smoother},
  Q. J. R. Meteorol. Soc., 140 (2014), pp.~1521--1535.

\bibitem{bougerol1993}
{\sc P.~Bougerol}, {\em Kalman filtering with random coefficients and
  contractions}, SIAM J. Control Optim., 31 (1993), pp.~942--959.

\bibitem{carrassi2008a}
{\sc A.~Carrassi, M.~Ghil, A.~Trevisan, and F.~Uboldi}, {\em Data assimilation
  as a nonlinear dynamical systems problem: Stability and convergence of the
  prediction-assimilation system}, Chaos, 18 (2008), p.~023112.

\bibitem{carrassi2008b}
{\sc A.~Carrassi, A.~Trevisan, L.~Descamps, O.~Talagrand, and F.~Uboldi}, {\em
  Controlling instabilities along a {3DVar} analysis cycle by assimilating in
  the unstable subspace: a comparison with the {EnKF}}, Nonlin. Processes
  Geophys., 15 (2008), pp.~503--521.

\bibitem{carrassi2009}
{\sc A.~Carrassi, S.~Vannitsem, D.~Zupanski, and M.~Zupanski}, {\em The maximum
  likelihood ensemble filter performances in chaotic systems}, Tellus A, 61
  (2009), pp.~587--600.

\bibitem{evensen2009}
{\sc G.~Evensen}, {\em {D}ata {A}ssimilation: {T}he {E}nsemble {K}alman
  {F}ilter}, Springer-Verlag Berlin Heildelberg, second~ed., 2009.

\bibitem{gajic1995}
{\sc Z.~Gaji{\'c} and M.~T.~J. Qureshi}, {\em Lyapunov matrix equation in
  system stability and control}, Academic Press, San Diego, California, 1995.

\bibitem{ginelli2013}
{\sc F.~Ginelli, H.~Chat{\'e}, R.~Livi, and A.~Politi}, {\em Covariant
  {L}yapunov vectors}, J. Phys. A: Math. Theor., 46 (2013), p.~254005.

\bibitem{gurumoorthy2017}
{\sc K.~S. Gurumoorthy, C.~Grudzien, A.~Apte, A.~Carrassi, and C.~K. R.~T.
  Jones}, {\em Rank deficiency of {K}alman error covariance matrices in linear
  time-varying system with deterministic evolution}, SIAM Journal on Control
  and Optimization, 0 (2017), pp.~0--0.
\newblock Accepted for publication, arXiv preprint arXiv:1503.05029.

\bibitem{jazwinski1970}
{\sc A.~H. Jazwinski}, {\em Stochastic Processes and Filtering Theory},
  Academic Press, New-York, 1970.

\bibitem{kalman1960}
{\sc R.~E. Kalman}, {\em A new approach to linear filtering and prediction
  problems}, Journal of Fluids Engineering, 82 (1960), pp.~35--45.

\bibitem{kuptsov2012}
{\sc P.~V. Kuptsov and U.~Parlitz}, {\em Theory and computation of covariant
  {L}yapunov vectors}, J. Nonlinear Sci., 22 (2012), pp.~727--762.

\bibitem{legras1996}
{\sc B.~Legras and R.~Vautard}, {\em A guide to lyapunov vectors}, in ECMWF
  Workshop on Predictability, Reading, United-Kingdom, 1996, ECMWF,
  pp.~135--146.

\bibitem{lorenz1998}
{\sc E.~N. Lorenz and K.~A. Emanuel}, {\em Optimal sites for supplementary
  weather observations: simulation with a small model}, J. Atmos. Sci., 55
  (1998), pp.~399--414.

\bibitem{ng2011}
{\sc G.-H.~C. Ng, D.~McLaughlin, D.~Entekhabi, and A.~Ahanin}, {\em The role of
  model dynamics in ensemble {K}alman filter performance for chaotic systems},
  Tellus A, 63 (2011), pp.~958--977.

\bibitem{palatella2013a}
{\sc L.~Palatella, A.~Carrassi, and A.~Trevisan}, {\em Lyapunov vectors and
  assimilation in the unstable subspace: theory and applications}, J. Phys. A:
  Math. Theor., 46 (2013), p.~254020.

\bibitem{palatella2015}
{\sc L.~Palatella and A.~Trevisan}, {\em Interaction of {L}yapunov vectors in
  the formulation of the nonlinear extension of the {K}alman filter}, Phys.
  Rev. E, 91 (2015), p.~042905.

\bibitem{palatella2013b}
{\sc L.~Palatella, A.~Trevisan, and S.~Rambaldi}, {\em Nonlinear stability of
  traffic models and the use of {L}yapunov vectors for estimating the traffic
  state}, Phys. Rev. E, 88 (2013), p.~022901.

\bibitem{parker1989}
{\sc T.~S. Parker and L.~O. Chua}, {\em Practical numerical algorithms for
  chaotic systems}, Springer-Verlag New York, 1989.

\bibitem{pham1998}
{\sc D.~T. Pham, J.~Verron, and M.~C. Roubaud}, {\em A singular evolutive
  extended {K}alman filter for data assimilation in oceanography}, J. Marine
  Syst., 16 (1998), pp.~323--340.

\bibitem{sakov2011}
{\sc P.~Sakov and L.~Bertino}, {\em Relation between two common localisation
  methods for the {E}n{KF}}, Comput. Geosci., 15 (2011), pp.~225--237.

\bibitem{simon2006}
{\sc D.~Simon}, {\em Optimal state estimation: {K}alman, {H} infinity, and
  nonlinear approaches}, John Wiley {\&} Sons, Inc., 2006.

\bibitem{trevisan2010}
{\sc A.~Trevisan, M.~D'Isidoro, and O.~Talagrand}, {\em Four-dimensional
  variational assimilation in the unstable subspace and the optimal subspace
  dimension}, Q. J. R. Meteorol. Soc., 136 (2010), pp.~487--496.

\bibitem{trevisan2011}
{\sc A.~Trevisan and L.~Palatella}, {\em On the {K}alman filter error
  covariance collapse into the unstable subspace}, Nonlin. Processes Geophys.,
  18 (2011), pp.~243--250.

\bibitem{trevisan1998}
{\sc A.~Trevisan and F.~Pancotti}, {\em Periodic orbits, {L}yapunov vectors,
  and singular vectors in the {L}orenz system}, J. Atmos. Sci., 55 (1998),
  pp.~390--398.

\bibitem{trevisan2004}
{\sc A.~Trevisan and F.~Uboldi}, {\em Assimilation of standard and targeted
  observations within the unstable subspace of the
  observation-analysis-forecast cycle}, J. Atmos. Sci., 61 (2004),
  pp.~103--113.

\bibitem{uboldi2006}
{\sc F.~Uboldi and A.~Trevisan}, {\em Detecting unstable structures and
  controlling error growth by assimilation of standard and adaptive
  observations in a primitive equation ocean model}, Nonlin. Processes
  Geophys., 16 (2006), pp.~67--81.

\bibitem{vannitsem2016}
{\sc S.~Vannitsem and V.~Lucarini}, {\em Statistical and dynamical properties
  of covariant {L}yapunov vectors in a coupled atmosphere-ocean
  model-multiscale effects, geometric degeneracy, and error dynamics}, J. Phys.
  A: Math. Theor., 49 (2016), p.~224001.

\bibitem{wojtowski2007}
{\sc M.~P. Wojtowski}, {\em Geometry of {K}alman fiters}, J.~Geom.~and Symmetry
  in Physics, 9 (2007), pp.~83--95.

\bibitem{wolfe2007}
{\sc C.~L. Wolfe and R.~M. Samelson}, {\em An efficient method for recovering
  lyapunov vectors from singular vectors}, Tellus A, 59 (2007), pp.~355--366.

\bibitem{zhang1999}
{\sc F.~Zhang}, {\em Matrix theory: basic results and techniques},
  Springer-Verlag New-York Inc., 1999.

\end{thebibliography}

\end{document}